\documentclass[preprint]{elsarticle}

\makeatletter
\def\ps@pprintTitle{%
 \let\@oddhead\@empty
 \let\@evenhead\@empty
 \def\@oddfoot{\reset@font\hfil\thepage\hfil}
 \let\@evenfoot\@oddfoot
}

\usepackage{graphics}
\usepackage{epsfig}
\usepackage{amsmath}
\usepackage{amsfonts}
\usepackage{cases}
\usepackage{amsthm}
\usepackage{color}
\usepackage{graphicx}
\usepackage{caption}
\usepackage{subcaption}
\usepackage{natbib}
\usepackage{placeins}

\catcode`@=11
\@addtoreset{equation}{section}
\renewcommand\theequation{\thesection.\@arabic\c@equation}
\catcode`@=12

\usepackage{geometry}
\usepackage{color}

\newcommand{\blue}[1]{{\color{blue} #1}}

\numberwithin{equation}{section}
\newtheorem{theorem}{\bf \blue{Theorem}}[section]
\newtheorem{definition}[theorem]{\bf \blue{Definition}}
\newtheorem{remark}[theorem]{\bf \blue{Remark}}
\newtheorem{example}[theorem]{\bf \blue{Example}}
\newtheorem{lemma}[theorem]{\bf \blue{Lemma}}

\usepackage{color}

\usepackage{graphicx}

\usepackage{color}

\usepackage{graphicx} 
\usepackage{indentfirst}

\usepackage[plainpages=false]{hyperref}
\hypersetup{colorlinks=true,citecolor=blue,filecolor=black,linkcolor=red,urlcolor=green}

\begin{document}
	
	\begin{frontmatter}

         \title{ \bf  Stochastic Mathematical Modelling Study for Understanding the Extinction, Persistence and Control of SARS-CoV-2 Virus at the Within-host Level
         \tnoteref{mytitlenote}}

          \author[nitt]{Bishal Chhetri}
        \ead{bishalc@iitk.ac.in}
    
        \author[nitt]{B.V. Ratish Kumar \corref{cor1}}
        \ead{bvrk@iitk.ac.in}

        \cortext[cor1]{Corresponding author}

		\address[nitt]{Department of Mathematics and Statistics, \\
  Indian Institute of Technology Kanpur, Kanpur - 208016, Uttar Pradesh, India}

\begin{abstract}
Stochastic differential equations characterized by uncertainty are effective in modelling virus dynamics and provide an alternative to traditional deterministic models. Epidemic models are inevitably subjected to the randomness within the system or the environmental noise. In this paper, we analyze the stochastic within-host compartment model for SARS-CoV-2 virus and explore its dynamics. We first examine the existence and positivity of the solution of the model using Ito’s formula and the establish the stochastic boundedness
and permanence of the model. Exponential stability of the infection-free equilibrium state is established. Numerical simulations are conducted to complement the theoretical results. Environmental noise is found to play a crucial role in the dynamics of the disease and can even lead to the extinction of the disease. The
model is also extended to a stochastic optimal control problem and the effectiveness of control measures, such as antiviral drugs and immunomodulators is investigated. 
\end{abstract}
		
		
\end{frontmatter}
	
\section{Introduction}\label{sec:intro}

The entire world faced a crisis due to the spread of the highly infectious SARS-CoV-2 virus. This virus was first detected in December 2019 in the city of Wuhan in China. On February 11, 2020, the World Health Organization (WHO) named the disease caused by the virus as COVID-19 \cite{tang2020hallmarks}. Millions of people have lost loved ones to this contagious virus, and deaths from its various strains are still being reported today. Symptoms of COVID-19 range from mild to severe, with the most common symptoms being coughing, phlegm, shortness of breath and fever \cite{weng2021pain}. These symptoms can vary from person to person, with older adults and those with other health conditions at higher risk of severe complications. SARS-CoV-2 is mainly transmitted by respiratory droplets.\\

In all infectious diseases, the epidemic process of disease transmission between hosts and the immunological process concerned with the virus-cell interactions play a very important role in understanding the dynamics of the disease. For decades, mathematical models have proven to provide useful information for understanding the behavior of infections across all levels.
Several mathematical models have been formulated and studied at the within-host (cellular) level as well as at the between-host levels \cite{al2016building, ben2015minimal, egonmwan2019mathematical,  Joshi2002, koutou2018mathematical, rachah2015mathematical, sene2020sir}. These studies have proven to be very effective in understanding the pattern of disease spread and developing effective optimal strategies to eliminate infections. The rapid spread of the SARS-CoV-2 virus posed a major challenge for health authorities around the world. The lack of understanding of the virus and its dynamics, as well as the lack of preparedness to deal with the disease, was a serious problem in the early phases of the disease. In order to understand and analyze the complex transmission pattern of the disease, a series of mathematical models were developed in different countries using the country-specific data and some understanding of the disease was gained. In \cite{cooper2020sir}, a deterministic compartment SIR model is developed  to investigate the effectiveness of the modeling approach on the pandemic. Studies describing the equilibrium solutions, stability analysis, and bifurcation analysis are discussed in \cite{biswas2020covid, chen2020mathematical, dashtbali2021compartmental, leontitsis2021seahir, ndairou2020mathematical, samui2020mathematical, sarkar2020modeling, zeb2020mathematical, khan2015stability, alqahtani2021mathematical, din2020mathematical}. Optimal control theory has also been used to study the effectiveness of different control strategies. Optimal control studies investigating the role and effectiveness of control measures such as treatment, vaccination, quarantine, isolation and screening are presented in
\cite{aronna2020model, dhaiban2021optimal, kkdjou2020optimal,libotte2020determination, ndondo2021analysis}. The role of partial lockdown to slow down the transmission of disease is discussed in \cite{aguiar2020modelling}.  \\

In mathematical biology, deterministic models, where the output of the model is  determined by the initial conditions, input parameters, and the mathematical equations that govern the system, are thoroughly researched and studied. Though deterministic models offers valuable insights by simulating real-world phenomena, they fall short in capturing the stochastic behaviours inherent in certain systems. Stochastic models are a special type of mathematical model that describe real-life phenomena characterized by uncertainty. Examples of systems influenced by such uncertainties include biological populations, traffic systems, stock markets, and complex computer networks. Mathematical epidemic models are inevitably subjected to the randomness within the system or the environmental noise, particularly when modelling biological phenomena such as viral dynamics \cite{tesfay2021dynamics}. Environmental noise can include factors such as, climate, habitats, medical quality, etc., which can alter the mortality rates, birth rates and transmission rates. Infectious diseases with the interactions between the host and the parasite is a complex system with components and processes acting at different levels of time, space and biological organization \cite{gutierrez2015within}.
The dynamics of these systems are directly or indirectly affected by many of these internal and external factors. Stochastic models add additional degree of realism into the modeling studies by taking into account the randomness factors that affect the system. This type of models represent the real world scenarios more closely and produce more realistic and meaningful outputs. Some of the epidemic modeling studies with stochastic effects are discussed in \cite{emvudu2016mathematical, liu2021stochastic,  ndii2017stochastic, zhang2020dynamics, allen2008introduction}. In \cite{sun2018asymptotic}, a delayed stochastic chemostat model is developed. Stochastic models for HIV virus and hepatitis virus are developed and studied in \cite{dalal2008stochastic, khan2018extinction} and the non-negative global solution and extinction conditions are derived. The conditions for the persistence and the extinction of a stochastic SIS model is studied in \cite{zhao2014threshold}. The authors in \cite{lahrouz2013extinction} have studied the existence of a stationary distribution for a stochastic SIRS model. A number of stochastic models are developed and studied on COVID-19 at the population level and some of these can be found in \cite{babaei2021mathematical,  hussain2022stochastic, tesfaye2021stochastic, srivastav2022deterministic}. \\

Within-host models focus on understanding viral dynamics and the interaction between viral particles and the host's immune response. These models have become crucial for decision making regarding therapeutic regimens and prophylactic interventions for various infectious and non-communicable diseases. Some of the within-host mathematical studies on COVID-19 can be found in \cite{hernandez2020host, almocera2021stability, elbaz2023viral}. Though many mathematical modeling studies are done on COVID-19 at the within-host level, most of them only focus on the deterministic dynamics ignoring the inherent randomness or the external noises. Inclusion of such randomness or noise into the models make them robust and reliable. So, in this study, we extend the work done in \cite{chhetri2021within} to stochastic framework by adding randomness into the model. The main focus of this present work is to investigate how environment noise affect within-host disease dynamics.  Initially, we introduce some basic concepts of stochastic differential equation and then formulate the stochastic model. Following this, we discuss the existence and uniqueness of the solution, boundedness and permanence, and the asymptotic stability. Subsequently, stochastic optimal control problem is formulated and optimal control solutions are derived using stochastic Maximum principle and Hamiltonian function \cite{bensoussan1988perturbation, el2021control}. The theoretical results obtained are validated with numerical illustrations.

\section{Basic Concepts of Stochastic Differential Equations}
\noindent
\begin{definition} 
\textnormal{Stochastic Processes}\\
\textnormal{A stochastic process is a parametrized collection of random variables $\{ X_t\}_{t\in T}$ defined on a probability space $(\Omega, F, P)$ and assuming values in $R^n$. Here, $\Omega$ is the sample space, and $F$ and $P$ denotes the $\sigma$ - algebra and the probability measure respectively}. 
\end{definition}

\noindent
\begin{definition}
\textnormal{Brownian Motion} \\
\textnormal{Brownian motion also known as Wiener processes is a continuous time stochastic process $W_t$ with the following properties.
\begin{enumerate}
    \item $W_0 = 0$ (at $t = 0)$ almost surely.
    \item The paths $t \to W(t)$ are continuous almost surely.
    \item It has independent and normally distributed increments. For $0\leq s \leq t$, $W_t - W_s$ is normally distributed with mean 0 and variance $t-s$.
\end{enumerate}
}
\end{definition}

\noindent
\begin{definition}
\textnormal{Stochastic Differential Equation (SDE)}\\
\textnormal{SDE takes the form:
\begin{equation}
    dX_t = b(t, X_t)dt + \sigma(t, X_t) d B_t \label{lm}
\end{equation}
Here, $b(t, X_t)$ represents the deterministic part of the system, $\sigma(t, X_t)$ represent the intensity of the noise, and $B_t$ is the Wiener process representing the white noise.}
\end{definition}

\noindent
\begin{definition}
\textnormal{1-D Ito's Formula}  \textnormal{\cite{oksendal2003stochastic, friedman1975stochastic}}\\
\textnormal{Let $X_t$ be Ito's process given by}
\begin{eqnarray*}
    d X_t = u dt + v d B_t
\end{eqnarray*}
\textnormal{Let $g(t, x) \in C^2([0, \infty)\times R)$. Then  $Y_t = g(t, X_t)$ is again an Ito's process and }
\begin{equation*}
    d Y_t = \frac{\partial g(t, X_t)}{\partial t}dt + \frac{\partial g(t, X_t)}{\partial x}dX_t +  \frac{1}{2} \frac{\partial g^2(t, X_t)}{\partial x^2}(dX_t)^2
\end{equation*}
\textnormal{where $(dX_t)^2 = dX_t . dX_t$ is computed according to the rules}
$dt.dt = dt.dB_t = dB_t.dt = 0 $,  $dB_t.d B_t = dt$\\
\noindent
\textnormal{Let $V$ denote family of all non-negative functions defined on $R^3 \times [0, \infty)$ such that they are twice continuously differentiable in $X$ and once in $t$. The differential operator $L$ for equation $(\ref{lm})$ is defined by}
\begin{equation*}
    L = \frac{\partial}{\partial t} + \sum_{i = 1}^d \frac{\partial}{\partial x_i}b_i + \frac{1}{2} \sum_{i, j = 1}^d [\sigma^T \sigma]\frac{\partial^2}{\partial x_i \partial x_j}
\end{equation*}
\noindent
\textnormal{If $L$ acts on V, then}
\begin{eqnarray*}
    L V(x, t) = V_t + V_x b + \frac{1}{2} trace \bigg[\sigma^T V_{xx} \sigma\bigg]
\end{eqnarray*}
\textnormal{where} $V_t = \frac{\partial V}{\partial t}$, \hspace{.3cm}
$V_x = \bigg(\frac{\partial V}{\partial x_1}, \hspace{.3cm}\frac{\partial V}{\partial x_2}, ...... \frac{\partial V}{\partial x_d}\bigg)$, \hspace{.3cm}
$V_{xx} = \frac{\partial^2 V}{\partial x_i \partial x_j}$. \\

\noindent
\textnormal{By using Ito's formula we get}
$$dV = L V(x, t) dt + V_x(x, t) \sigma(x, t) dB(t)$$
\end{definition}


\section{Stochastic Model Formulation}
In \cite{chhetri2021within}, a within-host mathematical model describing the interaction between susceptible cells (S), infected cells (I) and the viral particles (B) is developed based on the pathogenesis of the disease and the detailed study of the model is done. The model studied is given by the following system of equations $(\ref{1.1})-(\ref{1.3})$ with the parameter meaning in Table \ref{table11}. 

\begin{eqnarray}
   	\frac{dS}{dt}& =&  \omega \ - \beta SB  - \mu S  \label{1.1} \\
   	\frac{dI}{dt} &=& \beta SB \ -  {  p I }  \ - \mu I   \label{1.2}\\ 
   	\frac{dB}{dt} &=&  \alpha I   \ -  q B    \ -  \mu_{1} B \label{1.3}
   \end{eqnarray} 

The biologically feasible region of the system ($\ref{1.1})-(\ref{1.3})$ is given by, 
\begin{equation*}
\Omega = \bigg\{\bigg(S(t), I(t), B(t)\bigg) \in \mathbb{R}^{3}_{+} : 0 \leq S(t)  +  I(t) \leq N, \hspace{.2cm}0 \leq B(t)\leq B_{\text{max}}, \hspace{.12cm} \ t \geq 0 \bigg\}
\end{equation*}
where, $N$ is the maximum cell population in  a system and $B_{\text{max}}$ is the maximum bound on the viral particle.

\vspace{1cm}
\begin{table}[ht!]
     	\caption{Meanings of the Parameters}
     	\centering 
     	\begin{tabular}{|l|l|} 
     		\hline\hline
     		
     		\textbf{Parameters} &  \textbf{Biological Meaning} \\  
     		\hline\hline 
     	  $\omega$ & Natural birth rate of the susceptible cells \\
     	\hline\hline
     		$\alpha$ & Burst rate of the viral particles \\
     		\hline\hline
     		$\mu$ & Natural death rate of cells \\
     		\hline\hline
     	$\mu_1$ & Natural death rate of viral particles \\
     		\hline\hline
         $\beta$ & Transmission rate \\
         \hline\hline
         $p$ & Rate at which the infected cells are cleared \\ & by the individual immune response in the\\ & form of cytokines\\
          \hline\hline
          $q$ & Rate at which the virus particles are cleared \\ & by the individual immune response in the \\  & form of cytokines \\
          \hline\hline
         
     	\end{tabular} \label{table11}
     \end{table} \vspace{.25cm}

Mathematical epidemic models are influenced by environmental noise, with fluctuations in the environment significantly affecting model parameters and altering its dynamics \cite{tesfay2021dynamics, mahrouf2017dynamics}. Cell death can occur through various mechanisms, including apoptosis, necroptosis, and pyroptosis \cite{paolini2021cell}, all of which are influenced by different factors. The clearance of viral particles is also impacted by factors such as their binding and entry into cells, immune responses, and random the environmental noises \cite{mahrouf2017dynamics}. Given the complexity of biological phenomena and the external factors affecting the mortality rates, we incorporate noise into the within-host model $(\ref{1.1})-(\ref{1.3})$ to study the system's dynamics. Considering these factors into account we replace the parameters $\mu$ and $\mu_1$ of the model $(\ref{1.1})-(\ref{1.3})$ as follows.

$$\mu = \mu + \sigma_1 \frac{dW_1}{dt}$$
$$\mu_1 = \mu_1 + \sigma_2 \frac{dW_2}{dt}$$

Here $W_1(t)$,  $W_2(t)$ are mutually independent standard Brownian motions or Wiener process defined on a complete probability space $(\Omega, F, P)$ with a filtration $\{F_t \}_{t\geq 0}$.  We assume that the filtration $\{F_t\}_{t\geq 0}$ satisfies all the usual conditions (i.e., it is increasing and right continuous and $F_0$ contains all the $P$ null sets) \cite{oksendal2003stochastic}. The intensity of the noises are denoted by $\sigma_1$, and $\sigma_2$. Given that we assume the natural death rate for both uninfected and infected cells to be the same, the intensity of the noise, $\sigma_1$, and the Brownian motion $W_1$ are identical in the first two compartments. With this we assume that the biological factors or the external factors that affects the death of uninfected and infected cells are similar. 
With these assumption and using the updated definition of $\mu$ ad $\mu_1$, we get the following stochastic within-host model.

\begin{eqnarray}
   	{dS}& =&  \bigg(\omega \ - \beta SB  - \mu S\bigg)dt - \sigma_1 S dW_1(t)  \label{e1} \\
   	dI &=& \bigg(\beta SB - pI - \mu I \bigg) dt- \sigma_1 I dW_1(t)  \label{e2}\\ 
   	dB &=&  \bigg(\alpha I -qB - \mu_{1} B\bigg) dt - \sigma_2 BdW_2(t)\label{e3}
   \end{eqnarray}

\subsection{\textbf{Existence of Global Positive Solution}}
For every mathematical model, the first thing to prove is that the solution is non-negative. For a stochastic differential equation to have a unique global solution (meaning it doesn't blow up in finite time), the equation's coefficients usually need to meet two conditions: the linear growth condition and the local Lipschitz condition \cite{oksendal2003stochastic, yang2012ergodicity}. But because the interacting term $\beta S B$ is non linear,  the stochastic system $(\ref{e1})-(\ref{e3})$ fails to satisfy the linear growth condition. In almost all the interacting compartment models with noise linear growth condition fails \cite{yang2012ergodicity, witbooi2023stability, sengupta2018stochastic, zhang2020sufficient}. Therefore, in this case the solution of the stochastic system may explode at a finite time. To prove the existence of global positive solution Ito's formula is used and it is generally done using a contradiction approach.  Fort the deterministic model (\ref{1.1})-(\ref{1.3}) the coefficients of the model are shown to satisfy local Lipschitz condition and existence and uniqueness of the solutions are discussed in \cite{chhetri2021within, chhetri2022optimal}.\\

To simplify our analysis, we will make a change of notation here. Let $V_1 = S(t)$, $V_2 = I(t)$, and $V_3 = B(t)$. The viral particles is now denoted by $V_3(t)$. \\

Let $R_{+}^3 = \{V(t) = (V_1(t), V_2(t), V_3(t)) \in R^3: V_i > 0, 1\leq i \leq 3\} $. With $(V_1, V_2, V_3) = (S, I, B)$, system of equations $(\ref{e1})-(\ref{e3})$ becomes
\begin{eqnarray}
   	{dV_1}& =&  \bigg(\omega \ - \beta V_1V_3  - \mu V_1\bigg)dt - \sigma_1 V_1 dW_1(t)  \label{equ4} \\
   	dV_2 &=& \bigg(\beta V_1V_3- pV_2 - \mu V_2 \bigg) dt- \sigma_1 V_2 dW_1(t)  \label{equ5}\\ 
   dV_3 &=&  \bigg(\alpha V_2 -qV_3 - \mu_{1} V_3\bigg) dt - \sigma_2 V_3 dW_2(t)\label{equ6}
   \end{eqnarray} 

Before proving the main theorem, we will first prove a lemma.
\begin{lemma}
    The inequality $u \leq 2(u + 1 - log(u)) - (4-2 log 2)$ holds for all $u$ in R such that $u>0$ . \label{lemma}
\end{lemma}
\textbf{Proof:} To prove this lemma, let us define a function $f(u) = u + 2 - 2 log(u)$ for $u > 0$. Clearly, we find that $f(u)$ has a minimum at $u = 2$ and the minimum value is $4-2log2$. With simple calculations, we find that $u \leq 2(u + 1 - log(u)) - (4-2log2)$ holds $\forall u>0$.
We also note from the above lemma that for $u>0$, $u \leq 2(u + 1 - log(u))$.\\

\begin{theorem} 
For any initial value $V_0 \in R_+^3$, SDE $(\ref{equ4})-(\ref{equ6})$ has a unique solution $V(t)$ on $t\geq 0$ and solution will remain in $R_+^3$ with probability one for all $t\geq 0$ almost surely. \label{thm1}
\end{theorem}
\textbf{Proof:} It can be easily seen from equations $(\ref{equ4})-(\ref{equ6})$ that the coefficients are locally Lipschitz continuous. Therefore, for any given $V_0 \in R_+^3$, we have a unique solution for the SDE $(\ref{equ4})-(\ref{equ6})$ on $t \in [0, \tau_e)$, where $\tau_e$ is the explosion time ($\tau_e$ is a time at which $V(t)$ becomes infinite) \cite{mahrouf2017dynamics, din2021stochastic}. For the solution to be global, we need to show that $\tau_e = \infty$ almost surely. So in order to prove the theorem, all we need to show is $\tau_e = \infty$.\\

Let us take $k_0 \geq 0$ sufficiently large so that every component of $V_0$ lies within $[\frac{1}{k_0}, k_0]$. Now for each $k\geq k_0$,  let us define $\tau_k = inf \bigg \{t \in [0, \tau_e): V_i \notin \bigg(\frac{1}{k}, k\bigg) \text{for some i}, 1\leq i \leq 3 \bigg \}$.  $\tau_k$ is the stopping time at which any of $V_i$ leaves the interval $(\frac{1}{k}, k)$ for the first time. Let us set $inf(\phi) = \infty$ where $\phi$ is an empty set. Clearly we see that as $k$ becomes larger and larger $\tau_k$ increases. \\

Now, let $$\tau_\infty = \lim_{k \to \infty}\tau_k$$

By the definition of $\tau_k$, $\tau_\infty $ is a value between 0 and $\tau_e$. Therefore, in this case $\tau_\infty \leq \tau_e$. If we can show that $\tau_\infty = \infty$, then $\tau_e$ has to be equal to infinity and therefore we complete our proof. Let us assume that $\tau_\infty = \infty$ is false. In this case, we can find a pair of constants $0 < T $ and $\epsilon \in (0, 1)$ such that $P\{\tau_\infty \leq T\} \geq \epsilon$. Hence there is an integer $k_1 \geq k_0$ such that 

\begin{equation}
    P\{\tau_k \leq T\} \geq \epsilon \hspace{.5cm}\forall k \geq k_1 \label{probability}\\
\end{equation}

Define a $C^2$ function $\psi: R_+^3 \to R_+$ by
$$\psi(V) = \sum_{i=1}^{3} V_i + 1 - log V_i$$

We know that for $u> 0$, $u+1-log(u) > 0$. Therefore, by definition we see that $\psi(V)$ is a non-negative function. Using Ito formula on $\psi(V)$, we get the following expression.

\begin{equation*}
\begin{split}
{d \psi} & = \bigg[ (1 - \frac{1}{V_1})X + (1 - \frac{1}{V_2}) Y + (1 - \frac{1}{V_3}) Z + {\sigma_1^2}+ \frac{\sigma_2^2}{2}\bigg] dt + \sigma_1 (2 - V_1 - V_2)dW_1 + \sigma_2 (1-V_3)dW_2  \\[4pt]
& \leq \bigg[ \omega + 2\mu + \mu_1 + p + q + {\sigma_1^2} + \frac{\sigma_2^2}{2} + \beta V_3 + \alpha V_2 \bigg] dt + \sigma_1 (2 - V_1 - V_2)dW_1 + \sigma_2 (1-V_3)dW_2\end{split}
\end{equation*}
Here,
$$X = \omega - \beta V_1 V_3 -\mu V_1$$
$$Y = \beta V_1 V_3 - p V_2 - \mu V_2$$
$$Z = \alpha V_2 - q V_3 - \mu_1 V_3$$

 \noindent
Let 
$$C_1 = \omega + 2\mu + \mu_1 + p + q + {\sigma_1^2}+ \frac{\sigma_2^2}{2}$$
$$C_2 = 2(\alpha + \beta)$$
with $C_1$ and $C_2$ we have

\begin{equation}
d \psi \leq \bigg(C_1 + \alpha V_2 + \beta V_3\bigg)dt + \sigma_1 (2 - V_1 - V-2)dW_1 + \sigma_2 (1-V_3)dW_2  \label{eq3100}
\end{equation}

\noindent
From Lemma \ref{lemma} we know that $$V_i \leq 2 (V_i + 1 - log V_i)$$
Therefore,
$$\alpha V_2 \leq 2 \alpha (V_2 + 1 - log V_2)$$
$$\beta V_3 \leq 2 \beta (V_3 + 1 - log V_3)$$
\noindent
Using this we get
$$\alpha V_2 + \beta V_3 \leq 2\alpha(V_2 + 1 - log V_2) + 2 \beta (V_3 + 1 -log V_3)$$
\noindent
Calculating $C_2 \psi(V)$ we find that $$2\alpha(V_2 + 1 - log V_2) + 2 \beta (V_3 + 1 -log V_3) \leq C_2 \psi(V)$$

\noindent
Hence inequality (\ref{eq3100}) becomes
$${d \psi} \leq \bigg(C_1 + C_2 \psi(V)\bigg)dt + \sigma_1 (2 - V_1 - V_2)dW_1 + \sigma_2 (1-V_3)dW_2$$

\noindent
Let $C_3 =$ max$(C_1, \hspace{.2cm} C_2)$. With this the above inequality becomes
$${d \psi} \leq C_3 \bigg(1 + \psi(V)\bigg)dt + \sigma_1 (2 - V_1 - V_2)dW_1 + \sigma_2 (1-V_3)dW_2$$

\noindent
If $t_1 \leq T$ the integrating the above inequality we get\\

$\int_{0}^{\tau_k \wedge t_1} d \psi(V) \leq C_3 \int_{0}^{\tau_k \wedge t_1} \bigg(1 + \psi(V)\bigg)dt + \int_{0}^{\tau_k \wedge t_1} \sigma_1 (2 - V_1 - V_2)dW_1 + \int_{0}^{\tau_k \wedge t_1} \sigma_2 (1-V_3)dW_2 $\\

\vspace{.1cm}
\noindent
Here $\tau_k \wedge t_1 =$ min $(\tau_k, \hspace{.15cm}t_1)$
\begin{equation}
\begin{split}
\psi\bigg(V(\tau_k \wedge t_1)\bigg) &\leq \psi(V_0) + C_3 \int_{0}^{\tau_k \wedge t_1} \bigg(1 + \psi(V)\bigg)dt + \int_{0}^{\tau_k \wedge t_1} \sigma_1 (2 - V_1 - V_2)dW_1 + \int_{0}^{\tau_k \wedge t_1} \sigma_2 (1-V_3)dW_2 \\
&\leq \psi(V_0) + C_3 \int_{0}^{\tau_k \wedge t_1} \bigg(1 + \psi(V)\bigg)dt + \int_{0}^{\tau_k \wedge t_1} 2 \sigma_1dW_1 + \int_{0}^{\tau_k \wedge t_1} \sigma_2 dW_2 \\
& \leq \psi(V_0) + C_3 T +  C_3 \int_{0}^{\tau_k \wedge t_1} \psi(V)dt + 2 \sigma_1 W_1 + \sigma_2 W_2
\end{split}
\end{equation}
Now let us use expectation on both the sides.
$E\bigg[\psi\bigg(V(\tau_k \wedge t_1)\bigg)\bigg] \leq \psi(V_0) + C_3 T + C_3 E\bigg[\int_{0}^{\tau_k \wedge t_1} \psi(V)dt\bigg] $ with $E[W_i] = 0$.

\noindent
Using Gronwall's inequality \cite{gray2011stochastic} we get,
$$E\bigg[\psi\bigg(V(\tau_k \wedge t_1)\bigg)\bigg] \leq \bigg( \psi(V_0) + C_3 T\bigg)e^{C_3 T} = C_4$$.
This is the expected value of $\psi$ evaluated at the stopping time $\tau_k \wedge t_1$, i.e. minimum of the stopping time $\tau_k$ and a fixed time $t_1 \leq T$. \\

\noindent
Let $\Omega_k = \{\tau_k \leq T\}$ for $k\geq k_1$. $\Omega_k$ is the set of all sample paths $\omega_1$ for which the stopping time $\tau_k$
  is less than or equal to some fixed time $T$.  According to equation (\ref{probability}), we get $P(\Omega_k) \geq \epsilon$. This means, the probability that the solution $V_i$ for some $i, \hspace{.02cm} 1 \leq i\leq 3$  exits $(\frac{1}{k}, k)$ within $T$ is greater than a constant $\epsilon \in (0, 1)$.
Now for every $\omega_1 \in \Omega_k$, we can find $i  (1\leq i \leq 3)$ such that $V_i(\tau_k, \omega_1)$ equals either $k$ or $1/k$. In this case by definition $\psi(V(\tau_k, \omega_1))$ is not less than the smallest of $k + 1 - log k$ and $\frac{1}{k} + 1 + log k$.\\


\noindent
Therefore,
$$\psi\bigg(V(\tau_k, \omega_1)\bigg) \geq \bigg[ k + 1 - log k\bigg] \wedge \bigg[ \frac{1}{k} + 1 + log k\bigg]$$
$$C_4 \geq E\bigg[ \psi\bigg( V(\tau_k, \omega_1)\bigg)\bigg] = E\bigg[I_{\Omega_k} 
 \psi\bigg(V(\tau_k, \omega_1)\bigg)\bigg]\geq \epsilon \bigg[ k + 1 - log k\bigg] \wedge \bigg[ \frac{1}{k} + 1 + log k\bigg]$$
Here $I_{\Omega_k}$ is an indicator function of $\Omega_k$. As $k \to \infty$ both $( k + 1 -log k)$ and $(\frac{1}{k} + 1 + log k)$ tends to infinity. This leads to the contradiction $\infty > C_4 = \infty$. Therefore, the assumption $P(\tau_\infty \leq T) \geq \epsilon$ for $\epsilon \in (0, 1)$ is wrong. Hence, $\tau_\infty = \infty$ and this proves the theorem.\\

\begin{remark}
\textnormal{The above Theorem \ref{thm1} tells us that the solutions of stochastic system $(\ref{equ4})-(\ref{equ6})$ are positive and non-explosive. This may not be enough for a mathematical model to be meaningful. To show stochastic system $(\ref{equ4})-(\ref{equ6})$ is mathematically and biologically well-posed, it becomes essential to prove the property of stochastically ultimate boundedness. In the following we prove the boundedness and permanence in similar lines to that proved in \cite{zhang2023ultimate, cai2017stochastic}}.
\end{remark}

\subsection{\textbf{Stochastic Ultimate Boundedness}}
In this section, we will prove stochastic ultimate boundedness, which is another important property of the solution of stochastic differential system $(\ref{equ4})-(\ref{equ6})$. This property ensures that the state variable (population count in our case) remains within a certain bounds over time. This property is very crucial in many applications, especially in the areas of control theory and stability analysis, to ensure that the system does not explode over time.

\begin{definition}\textnormal{\cite{bahar2004stochastic, li2009population}}
\textnormal{The solution $X(t)$ of system $(\ref{equ4})-(\ref{equ6})$ is said to be stochastically ultimately bounded, if for any $\epsilon \in (0, 1)$, there is a positive constant $k = k(\epsilon)$, such
that for any initial values, the solution $X(t) = \bigg(V_1(t), V_2(t), V_3(t)\bigg)$ of the system $(\ref{equ4})-(\ref{equ6})$ has the property that $${lim_{t \rightarrow{\infty}} sup\hspace{.1cm}P}\bigg\{|X(t)| > k\bigg\} < \epsilon $$}
\label{def1}
\end{definition}

\begin{theorem}
    The solutions of the stochastic differential equations $(\ref{equ4})-(\ref{equ6})$ are stochastically ultimately bounded for any initial value $X_0 \in R_+^3$ .
\end{theorem}
\noindent
\textbf{Proof.}
Let us define a fuinction
$$\phi(X) = V_1 + V_2 + V_3$$ and $$X = (V_1, V_2, V_3)$$
Clearly, $\phi(X) \geq 0$ and $\phi(X) \rightarrow \infty$ as $|X| \rightarrow \infty$.
Using Ito's formula we get
\begin{equation}
    d \phi(t) = L \phi(t) dt - \sigma_1 V_1 dW_1 - \sigma_1 V_2 dW_2 - \sigma_2 V_3 dW_3
\end{equation}
where
\begin{eqnarray*}
    L \phi(t) = \omega + (\alpha - p) V_2 - (q + \mu_1)V_3 -\mu(V_1+V_2)
\end{eqnarray*}
Now similarly we have
\begin{equation}
\begin{split}
     d(e^{\mu t} \phi) &=L e^{\mu t} \phi dt - e^{\mu t}\bigg\{ \sigma_1 V_1 dW_1 + \sigma_1 V_2 dW_2 + \sigma_2 V_3 dW_3\bigg\}\\
     &= e^{\mu t} \bigg \{ \omega -(p - \alpha) V_2 - q V_3 -(\mu_1 - \mu)V_3   \bigg\}dt - e^{\mu t} \bigg\{ \sigma_1 V_1 dW_1 + \sigma_1 V_2 dW_2 + \sigma_2 V_3 dW_3\bigg\}  \label{eq310}
\end{split}
\end{equation}
Integrating the above equation $(\ref{eq310})$ from $0$ to $t$ we get
\begin{equation}
\begin{split}
e^{\mu t} \phi(t) &= \phi(0) + e^{\mu t} \omega - \omega - \int_{0}^ t e^{\mu s}(p-\alpha) V_2(s)ds - q \int_{0}^t e^{\mu s} V_3(s)ds - (\mu_1 - \mu) \int_{0}^t e^{\mu s} V_3(s)ds \\
& - \int_{0}^t \sigma_1 e^{\mu s} V_1(s) dW_1  - \int_{0}^t \sigma_1 e^{\mu s} V_2(s) dW_1 - \int_{0}^t \sigma_2 e^{\mu s} V_3(s) dW_2 
    \end{split}
\end{equation}
Taking expectation on both the sides
\begin{equation}
    \begin{split}
        e^{\mu t} E[\phi(t)] &= \phi(0) + \omega(e^{\mu t}-1) - E \int_{0}^ t e^{\mu s}(p-\alpha) V_2(s)ds - q E \int_{0}^t e^{\mu s} V_3(s)ds - (\mu_1 - \mu) E \int_{0}^t e^{\mu s} V_3(s)ds \\
        & - E \int_{0}^t \sigma_1 e^{\mu s} V_1(s) dW_1  - E \int_{0}^t \sigma_1 e^{\mu s} V_2(s) dW_1 - E \int_{0}^t \sigma_2 e^{\mu s} V_3(s) dW_2 \\
        &\leq \phi(0) + \omega(e^{\mu t}-1) - E \int_{0}^ t e^{\mu s}(p-\alpha) V_2(s)ds - q E \int_{0}^t e^{\mu s} V_3(s)ds - (\mu_1 - \mu) E \int_{0}^t e^{\mu s} V_3(s)ds
    \end{split} \label{312}
\end{equation}
The equation above (\ref{312}) is derived from the property that the increments $dW_t$ of a Wiener process have a zero mean by definition. As a result, the integral of any function with respect to $dW_t$ will also have a zero mean \cite{oksendal2013stochastic, klebaner2012introduction}.\\

\noindent
\noindent
Since the solution $V_1(t), V_2(t)$ and $V_3(t)$ of system $(\ref{equ4})-(\ref{equ6})$ will not tend to infinity in finite time, so there must be a constant number $K_1$ such that 
$$e^{\mu t} E [\phi(t)] \leq \phi(0)+\omega(e^{\mu t}-1) + K_1$$
as $t \rightarrow \infty$ we get
$$E[\phi(t)] \leq \omega$$ and
$$\text{lim}_{t \rightarrow \infty} \hspace{.15cm}\text{sup} \hspace{.08cm}E [\phi(t)] \leq \omega$$

\noindent
Since $X(t) = (V_1(t), V_2(t), V_3(t))$
\begin{equation*}
    \begin{split}
        |X(t)|^2 &= V_1^2 + V_2^2 + V_3^2 \leq 3 \hspace{.1cm}\text{max}\{V_1^2, V_2^2, V_3^2\}
    \end{split}
\end{equation*}
Therefore
\begin{equation*}
    \begin{split}
        |X(t)| &\leq \sqrt{3} \hspace{.3cm} \text{max}\{V_1, V_2, V_3\}\\
        & \leq \sqrt{3} \hspace{.1cm} (V_1 + V_2 + V_3) = \sqrt{3} \hspace{.1cm} \phi(t)
    \end{split}
\end{equation*}
Now
$$\text{lim}_{t \rightarrow \infty} \hspace{.15cm}\text{sup} \hspace{.2cm}E [\phi(t)] \leq \sqrt{3} \hspace{.3cm}\text{lim}_{t \rightarrow \infty} \hspace{.15cm}\text{sup} \hspace{.3cm}E[\phi(t)] \leq \sqrt{3} \hspace{.3cm}\omega $$
Therefore
$$\text{lim}_{t \rightarrow \infty} \hspace{.15cm}\text{sup} E[X(t)] \leq \sqrt{3} \hspace{.3cm}\omega $$
For any $\epsilon > 0$, let us choose $k(\epsilon) = \frac{\sqrt{3} \hspace{.1cm}\omega}{\epsilon}$\\

\noindent
Using Chebyshev's inequality \cite{zhang2023ultimate, cai2017stochastic} we get
$$P \bigg(|X| \geq k(\epsilon)\bigg) \leq \frac{E[|X(t)|]}{k(\epsilon)}$$
In other words
$$\text{lim}_{t \rightarrow \infty} \hspace{.15cm}\text{sup} \hspace{.15cm} P \bigg(|X| \geq k(\epsilon)\bigg) \leq \frac{\sqrt{3}\hspace{.1cm}\omega}{\frac{\sqrt{3}\hspace{.1cm}\omega}{\epsilon}}= \epsilon$$
Hence, by Definition \ref{def1} the stochastic system $(\ref{equ4})-(\ref{equ6})$ is
stochastically ultimately bounded. \\

\subsection{\textbf{Stochastic Permanence}}
In this section we show that the system $(\ref{equ4})-(\ref{equ6})$ is stochastically permanent. This property ensures that the variables of the system do not approach zero or infinity over time. The concept of system permanence is a critical technical issue, relevant to various types of systems such as social, medical, biological, population, mechanical, and electrical systems. \\

\begin{definition}\textnormal{\cite{cai2017stochastic, jiang2008global}}
\textnormal{The model $(\ref{equ4})-(\ref{equ6})$ is said to be stochastically permanent, if for any $\epsilon \in (0, 1)$, there exists a pair of positive constants $\rho = \rho(\epsilon)$ and $\xi = \xi(\epsilon)$ such
that for any initial values $X_0 \in R_{+}^3$, the solution $X(t) = \bigg(V_1(t), V_2(t), V_3(t)\bigg)$ of the system $(\ref{equ4})-(\ref{equ6})$ has the property that} 
$${lim_{t \rightarrow{\infty}}\hspace{.1cm}inf\hspace{.1cm}P}\bigg\{|X(t)| \leq \rho\bigg\} \geq 1- \epsilon $$ and
$${lim_{t \rightarrow{\infty}}\hspace{.1cm}inf\hspace{.1cm}P}\bigg\{|X(t)| \geq \xi \bigg\} \geq 1- \epsilon $$
\label{def2}
\end{definition}

\begin{theorem}
    The solutions of the stochastic differential equations $(\ref{equ4})-(\ref{equ6})$ are stochastically permanent for any initial value $X_0 \in R_+^3$ .
\end{theorem}
\noindent
\textbf{Proof}
Let us define $N(t) = V_1(t) + V_2(t) + V_3(t)$ and
$\Phi(t) = N + \frac{1}{N}$.\\

\noindent
Using Ito's formula we get
\begin{equation}
    d(e^{\mu t} \Phi) = L(e^{\mu t} \Phi(t)) dt - e^{\mu t} \{ \sigma_1 V_1 dW_1 + \sigma_1 V_2 dW_1 + \sigma_1 V_3 dW_2   \} \label{eq313}
\end{equation}
where
\begin{equation*}
\begin{split}
    L(e^{\mu t} \Phi(t)) &= \mu e^{\mu t} \Phi + e^{\mu t}\bigg[ (\omega - \mu(V_1 + V_2) + (\alpha - p)V_2 - (q+\mu_1)V_3) - \frac{1}{N^2}(\omega - \mu(V_1 + V_2) \\
    &+ (\alpha - p)V_2 - (q+\mu_1)V_3) \bigg] + \frac{e^{\mu t}}{N^3}(\sigma_1^2 V_1^2 + \sigma_1^2 V_2^2 + \sigma_1^2 V_3^2)
\end{split}
\end{equation*}
\noindent
 Integrating equation (\ref{eq313}) from $0$ to $t$ and then taking expectation on both the side we get
 \begin{equation}
     E[e^{\mu t} \Phi(t)] \leq \Phi(0) + \frac{\omega}{\mu}(e^{\mu t} - 1) + K
 \end{equation}
 where $K$ is a constant. Note that since the solution of system  $(\ref{equ4})-(\ref{equ6})$ will not tend to infinity in finite
time and remain positive, we get a constant $K$ after integrating and taking the expectation on the remaining terms of equation (\ref{eq313}). 
\begin{equation}
    e^{\mu t} E[ \Phi(t)] \leq \Phi(0) + \frac{\omega}{\mu}(e^{\mu t} - 1) + K
 \end{equation}
as $t$ approaches $\infty$ we get
$$E[ \Phi(t)]  \leq\frac{\omega}{\mu} = H$$
Now let $\epsilon > 0$ and choose constant $\rho$ sufficiently large such that $\frac{H}{\rho} < 1$. By Chebyshev's inequality
$$P(N + \frac{1}{N} > \rho) \leq \frac{E[N + \frac{1}{N}]}{\rho} \leq \frac{H}{\rho} = \epsilon$$
Therefore
$$P(N + \frac{1}{N} > \rho) \leq \epsilon$$
\noindent
This implies that
$$P(N + \frac{1}{N} \leq \rho) \geq 1 - \epsilon$$
Or,
\begin{equation}
   1-\epsilon \leq P(N + \frac{1}{N} \leq \rho) \leq P(\frac{1}{\rho} \leq N \leq \rho) \label{ineq} 
\end{equation}
\noindent
Now, since $N = V_1 + V_2 + V_3$ and $X = (V_1, V_2, V_3)$
$$|X|^2 = V_1^2 + V_2^2 + V_3^2$$
$$N^2 = (V_1 + V_2 + V_3)^2 \leq 3(V_1^2 + V_2^2 + V_3^2) = 3 |X|^2$$
Therefore from the above we get
$$N^2 \leq 3 |X|^2 \leq 3 N^2$$
Or,
$$\frac{N}{\sqrt{3}} \leq |X| \leq \sqrt{3} N$$
From inequality (\ref{ineq}) we have
\begin{equation*}
    P\bigg(\frac{1}{\sqrt{3} \rho } \leq \frac{N}{\sqrt{3}} \leq |X| \leq N \leq \rho\bigg) \geq 1 - \epsilon
\end{equation*}
Therefore by Definition \ref{def2} the stochastic differential equations $(\ref{equ4})-(\ref{equ6})$ are stochastically permanent. The proof is complete.

\subsection{{\textbf{Asymptotic Analysis}}}\vspace{.25cm}	
 The basic reproduction number for the deterministic part of the system $(\ref{equ4})-(\ref{equ6})$$(\sigma_i = 0)$ is given by 

\begin{equation}
{ R_{0}}= {\frac{\beta \alpha \omega}{\mu (p+\mu) (q+\mu_{1})}} \label{R0}
\end{equation}

Deterministic system is shown to have two equilibrium points namely, the infection free equilibrium $E_{0}=\bigg(\frac{\omega}{\mu},0,0 \bigg)$ and the infected equilibrium $E_{1}=(V_1^*, V_2^*, V_3^*)$  where,

$$\hspace{.3cm}V_1^*=\frac{(p+\mu)(q+\mu_{1})}{\alpha \beta}$$
$$\hspace{.8cm}V_2^*=\frac{\mu(\mu_{1}+ q) \bigg(R_{0}-1\bigg)}{\alpha \beta}$$
$$\hspace{-.4cm}V_3^*=\frac{\mu \bigg(R_{0}-1\bigg)}{\beta}$$

In this case the infected equilibrium exists only if $R_{0}>1$. The detailed analysis for the stability of these equilibrium points are discussed in \cite{chhetri2021within}. In this study, we focus on the exponential stability of the infected cell population $(V_2 (t))$ and the viral particles $(V_3(t))$ for the SDE  $(\ref{equ4})-(\ref{equ6})$. \\

\begin{theorem}
    The infected cell and viral load population $V_2(t)$ and $V_3(t)$ of the system $(\ref{equ4})-(\ref{equ6})$ are exponentially stable almost surely in the sense that $V_2(t)$ and $V_3(t)$ will tend to their equilibrium
value 0 with probability 1 if the following two conditions are satisfied.\\
A). $\bigg(\alpha + \frac{\beta \omega }{\mu}\bigg) - (p+q+\mu+\mu_1) < \frac{\sigma_1^2 +\sigma_2^2}{2}$ \\
B). $\bigg(\alpha + \frac{ \beta \omega}{\mu} - (p +\mu + q+\mu_1)\bigg) < [2(\alpha -p-\mu)-\sigma_1^2][2\frac{\beta \omega}{\mu} - 2(q+\mu_1) - \sigma_2^2]$
 \label{thm38}
\end{theorem}	

\textbf{Proof} By definition \cite{gray2011stochastic}, the trivial solution of SDE $(\ref{lm})$ is said to be exponentially stable almost surely if
$$ lim_{k \to \infty}sup\frac{1}{t}log(x(t;t_0, x_0)) < 0$$. 
We will use this definition and show that the solution $V_1(t)$ and $V_2(t)$ of the system (\ref{equ4})-(\ref{equ6}) are exponentially stable. We assume that the parameters of the model are chosen such that the conditions A) and B) are satisfied. From the SDE $(\ref{equ4})-(\ref{equ6})$ we have
\begin{eqnarray}
   	dV_2 &=& \bigg(\beta V_1V_3- pV_2 - \mu V_2 \bigg) dt- \sigma_1 V_2 dW_1(t)  \label{e8}\\ 
   dV_3 &=&  \bigg(\alpha V_2 -qV_3 - \mu_{1} V_3\bigg) dt - \sigma_2 V_3 dW_2(t)\label{e9}
   \end{eqnarray}  

Now
$$d(V_2 + V_3) = \bigg[ \beta V_1 V_2 - p V_2 - \mu V_2 + \alpha V_2 - q V_3 -\mu_1 V_3 \bigg]dt - \sigma_1 V_2 dW_1 - \sigma_2 V_3 dW_2$$
Let $V=(V_2, V_3)$ and let us define a function $\Phi(V) = log(V_2 + V_3)$ for $V_2, V_3 \in (0, \infty)$.\\

\noindent
Using Ito's formula on $\Phi(V)$ we get the following.
\begin{equation}
\begin{split}
   d\Phi(V) &= \frac{1}{2(V_2 + V_3)^2} \bigg[ 2(V_2 + V_3)\bigg(\beta V_1 V_2 - pV_2 - \mu V_2 + \alpha V_2 - qV_3 - \mu_1 V_3 \bigg) - \sigma_1^2 V_2^2 - \sigma_2^2 V_3^2\bigg] dt \\ & - \frac{\sigma_1 V_2}{V_2 + V_3}dW_1 - \frac{\sigma_2 V_3}{V_2 + V_3}dW_2  \label{e12} \\
   & \leq \frac{1}{2(V_2 + V_3)^2} \bigg[ 2(V_2 + V_3)\bigg(\frac{\beta \omega}{\mu} V_2 - pV_2 - \mu V_2 + \alpha V_2 - qV_3 - \mu_1 V_3 \bigg) - \sigma_1^2 V_2^2 - \sigma_2^2 V_3^2\bigg] dt \\ & - \frac{\sigma_1 V_2}{V_2 + V_3}dW_1 - \frac{\sigma_2 V_3}{V_2 + V_3}dW_2
   \end{split}
\end{equation}
The above inequality results from the fact that $V_1(t)$ is bounded above by $\frac{\omega}{\mu}$. \\ The expression $2(V_2 + V_3)\bigg(\frac{\beta \omega}{\mu} V_2 - pV_2 - \mu V_2 + \alpha V_2 - qV_3 - \mu_1 V_3 \bigg) - \sigma_1^2 V_2^2 - \sigma_2^2 V_3^2$ can be simplified and written in terms of a matrix as,
$VAV^T$ where $V=(V_2, V_3)$ and 
\begin{equation*}
A = 
\begin{pmatrix}
2(\alpha - p - \mu) - \sigma_1^2 & \alpha + \frac{\beta \omega}{\mu} - p-\mu - q- \mu_1 \\
\alpha + \frac{\beta \omega}{\mu} - p-\mu - q- \mu_1 & 2 \frac{\beta \omega}{\mu} - 2(q+\mu_1) - \sigma_2^2 \\
\end{pmatrix}
\end{equation*} 

\noindent
Matrix $A$ is symmetric and with conditions A) and B) both the eigenvalues $\lambda_1$ and $\lambda_2$ of matrix $A$ are real and negative. \\

\noindent
Let $\lambda_{max} = $max$(\lambda_1, \lambda_2)$. $A$ is a negative definite symmetric matrix with $\lambda_{max}$ as the largest eigenvalue. Therefore we have the following inequality.
$$V^T A V \leq \lambda_{max} (V_2^2 + V_3^2) = - |\lambda_{max}| (V_2^2 + V_3^2)$$
\noindent
 Now equation $(\ref{e12})$ becomes
\begin{equation}
\begin{split}
   d\Phi(V) \leq \frac{-|\lambda_{max}|}{2(V_2 + V_3)^2}  (V_2^2 + V_3^2 )dt - \frac{\sigma_1 V_2}{V_2 + V_3}dW_1 - \frac{\sigma_2 V_3}{V_2 + V_3}dW_2  \label{e13} 
   \end{split}
\end{equation}
\noindent
Since $V_2 V_3 \leq \frac{V_2^2 + V_3^2}{2}$, we have $(V_2 + V_3)^2 \leq 2(V_2^2 + V_3^2)$ and -$(V_2^2 + V_3^2)\leq \frac{-(V_2 + V_3)^2}{2}$. Using this inequality we have 

\begin{equation}
\begin{split}
   d\Phi(V) \leq \frac{-|\lambda_{max}|}{4}dt - \frac{\sigma_1 V_2}{V_2 + V_3}dW_1 - \frac{\sigma_2 V_3}{V_2 + V_3}dW_2  \label{e14} 
   \end{split}
\end{equation}

\noindent
Since $\Phi(V) = log(V_2 + V_3)$. From $(\ref{e14})$ we get
$$d(log(V_2 + V_3)) \leq \frac{-|\lambda_{max}|}{4}dt - \frac{\sigma_1 V_2}{V_2 + V_3}dW_1 - \frac{\sigma_2 V_3}{V_2 + V_3}dW_2 $$
Integrating both the side,
$$log(V_2 + V_3) \leq \frac{-|\lambda_{max}|}{4}t - \int \frac{\sigma_1 V_2}{V_2 + V_3}dW_1 - \int \frac{\sigma_2 V_3}{V_2 + V_3}dW_2$$
Dividing by $t$ we have,
$$\frac{1}{t}log(V_2 + V_3) \leq \frac{-|\lambda_{max}|}{4} - \frac{1}{t}\int \frac{\sigma_1 V_2}{V_2 + V_3}dW_1 - \frac{1}{t} \int \frac{\sigma_2 V_3}{V_2 + V_3}dW_2$$
\noindent
Using large number theorem \cite{mao2007stochastic} we have
$$lim_{t\to \infty} sup \frac{1}{t}\int \frac{\sigma_1 V_2}{V_2 + V_3}dW_1 = lim_{t\to \infty} sup \frac{1}{t} \int \frac{\sigma_2 V_3}{V_2 + V_3}dW_2 = 0$$
Now
$$ lim_{t \to \infty}sup\frac{1}{t}log(V_2 + V_3) \leq \frac{-|\lambda_{max}|}{4} < 0$$
 Hence we conclude that as $t \to \infty$, $V_2(t)$ and $V_3(t)$ both tend to their disease-free equilibrium point 0 almost surely.

\begin{remark}
\textnormal{The above Theorem \ref{thm38} tells us that the infected cell population and the viral load population decays exponentially and goes to zero almost surely with appropriate choice of the parameter values. We infer from this that larger noises can lead to the extinction of the disease. }
\end{remark}

\section{\textbf{Numerical Illustrations}}
  In the above Theorem \ref{thm38}, we showed that $V_2(t)$ and $V_3(t)$ tend to the disease free state $(0, 0)$ exponentially with the two conditions (A, B) being satisfied. In this section, we support these results by simulations. We used MATLAB for the simulations and the results were verified by running the simulation repeatedly and extensively. Euler-Maruyama method is used for solving the stochastic differential equation. This method is a simple and effective numerical method for approximating the solutions of stochastic differential equations (SDEs) \cite{bayram2018numerical, higham2001algorithmic}. For the stochastic differential equation
$$dX(t)=f(X(t),t)dt+g(X(t),t)dW(t)$$
Euler-Maruyama method iterates over small time steps $\delta t$ using the following iterative formula.
$$X_{n+1} = X_n + f(X_n, t_n)\Delta t + g(X_n, t_n)\Delta W_n$$
$X_n$ is an approximation of $X(t_n)$ at time $t_n$ where $t_n = t_0 + n\Delta t$. $\Delta W_n = W(t_{n +1})- W(t_n)$ are independent  normal random variables with mean 0 and variance $\Delta t$, $\Delta t = \frac{(T - t_0)}{N}$ for $N$ time steps. The parameter values used in the simulations are all taken from the published works and the details are given in the following table. 
  
  \begin{table}[ht!]
     	\caption{parameter Values}
     	\centering 
     	\begin{tabular}{|l|l|l|} 
     		\hline\hline
     		
     		\textbf{Symbols} &  \textbf{Values}& \textbf{Source} \\  
     		\hline\hline 
     	  $\omega$ & $10$ & \cite{chhetri2021within} \\
     	\hline\hline
      $\beta$ & $0.005$ & \cite{chhetri2021within}\\
     	\hline\hline
     	 $\mu$ & $0.1 $& \cite{chhetri2021within}\\
     	\hline\hline
     		 $\mu_1$ & $0.6$ &\cite{chhetri2021within} \\
     		\hline\hline
     		 $\alpha$ & $0.24 $& \cite{li2020within} \\
     		 \hline\hline
     		 $p$ & $0.795$&  \cite{chhetri2021within} \\
     		 \hline\hline
     	 $q$ & $0.28$&  \cite{ghosh2021within} \\
     		 \hline\hline
     			
\end{tabular} \label{table2}
     \end{table} 

\begin{figure}[hbt!]
\begin{center}
\includegraphics[width=3.8in, height=3.2in, angle=0]{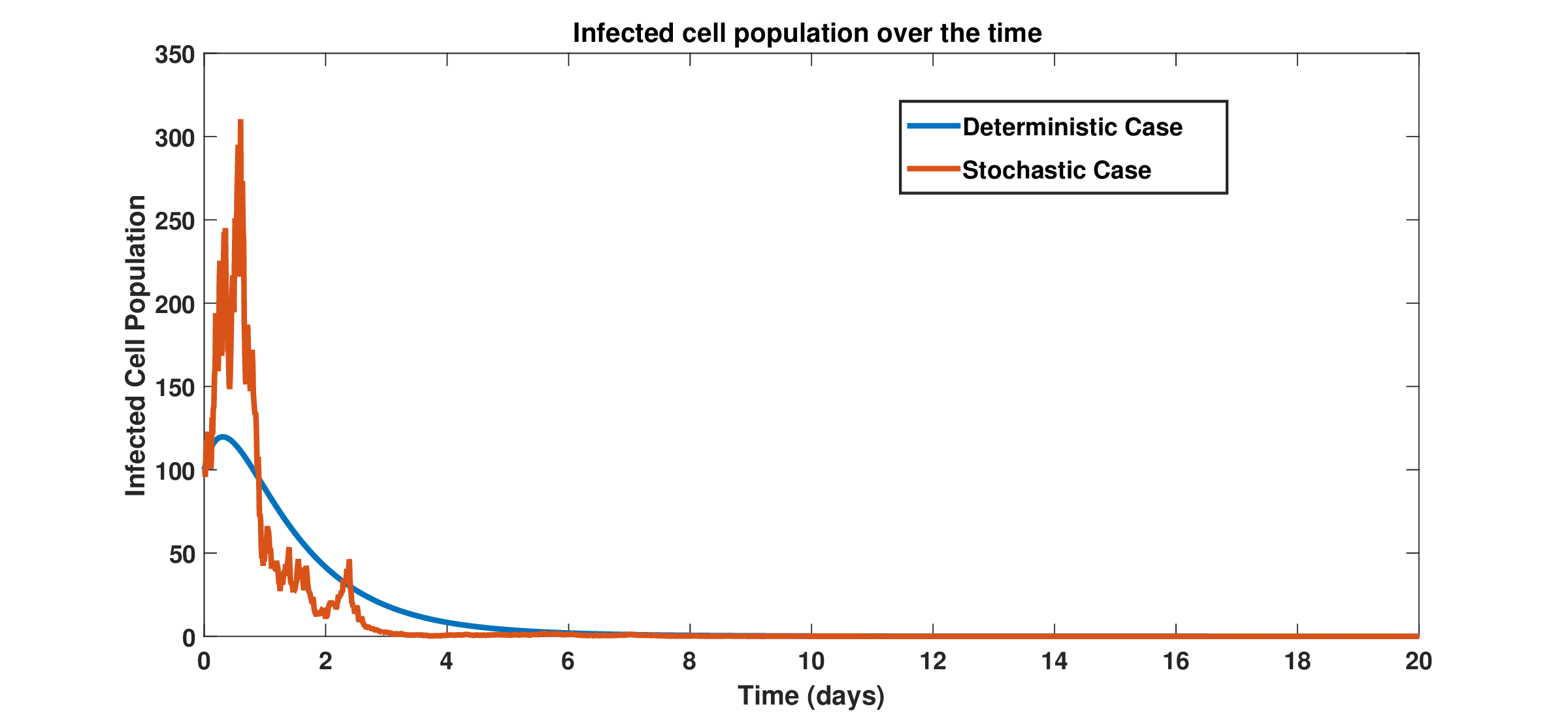}
\hspace{-.4cm}
\includegraphics[width=3.8in, height=3.2in, angle=0]{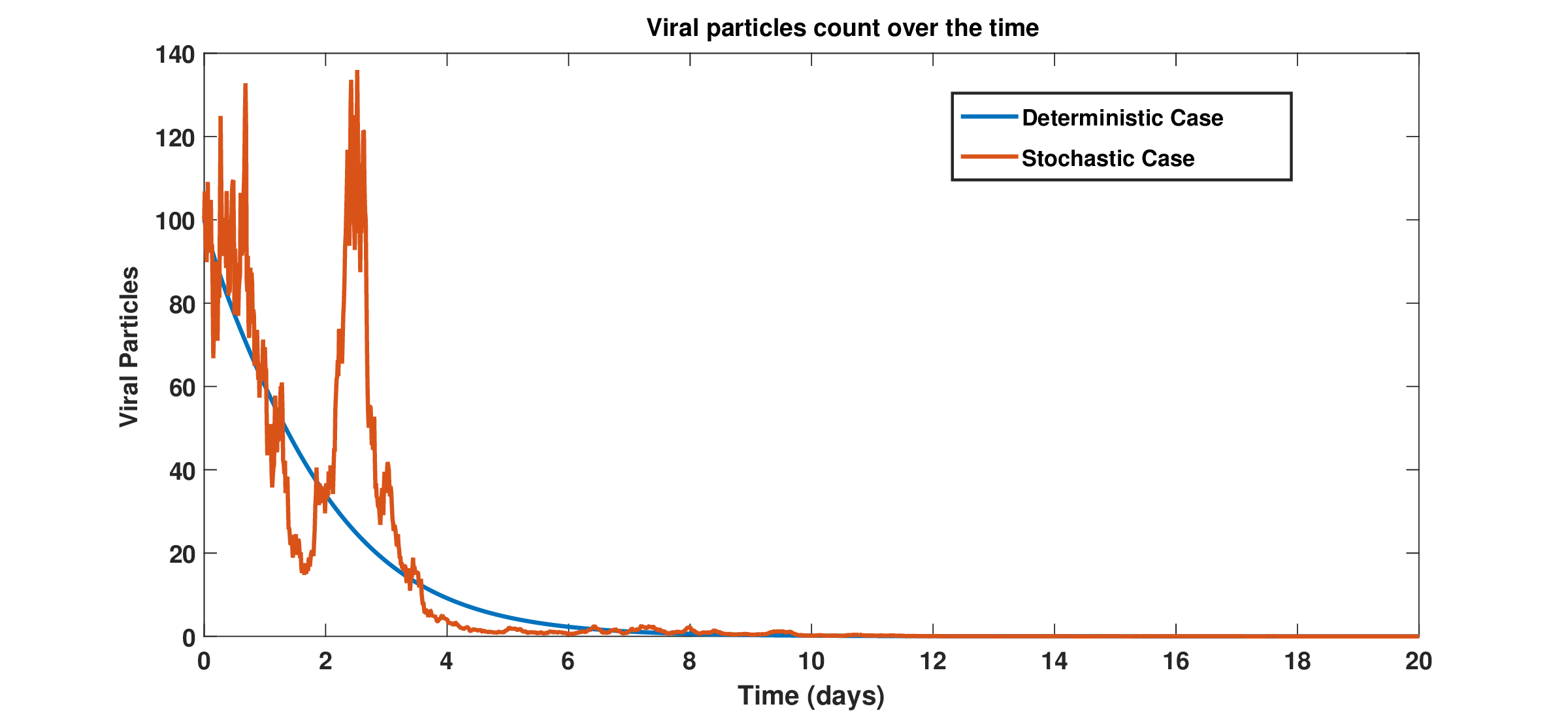}
\vspace{5mm}
\caption{Figure depicting the infected cell population and the viral load for the deterministic and stochastic cases. The solution of the system $(\ref{equ4})-(\ref{equ6})$ tends to the infection-free state $(I^*, B^*) = (0, 0)$ with $\sigma_1 = \sigma_2 = 0.1$ and other parameters values from Table \ref{table2} satisfying the conditions of Theorem \ref{thm38}.}
\label{fig1}
\end{center}
\end{figure}

\begin{figure}[hbt!]
\begin{center}
\includegraphics[width=3.8in, height=3.3in, angle=0]{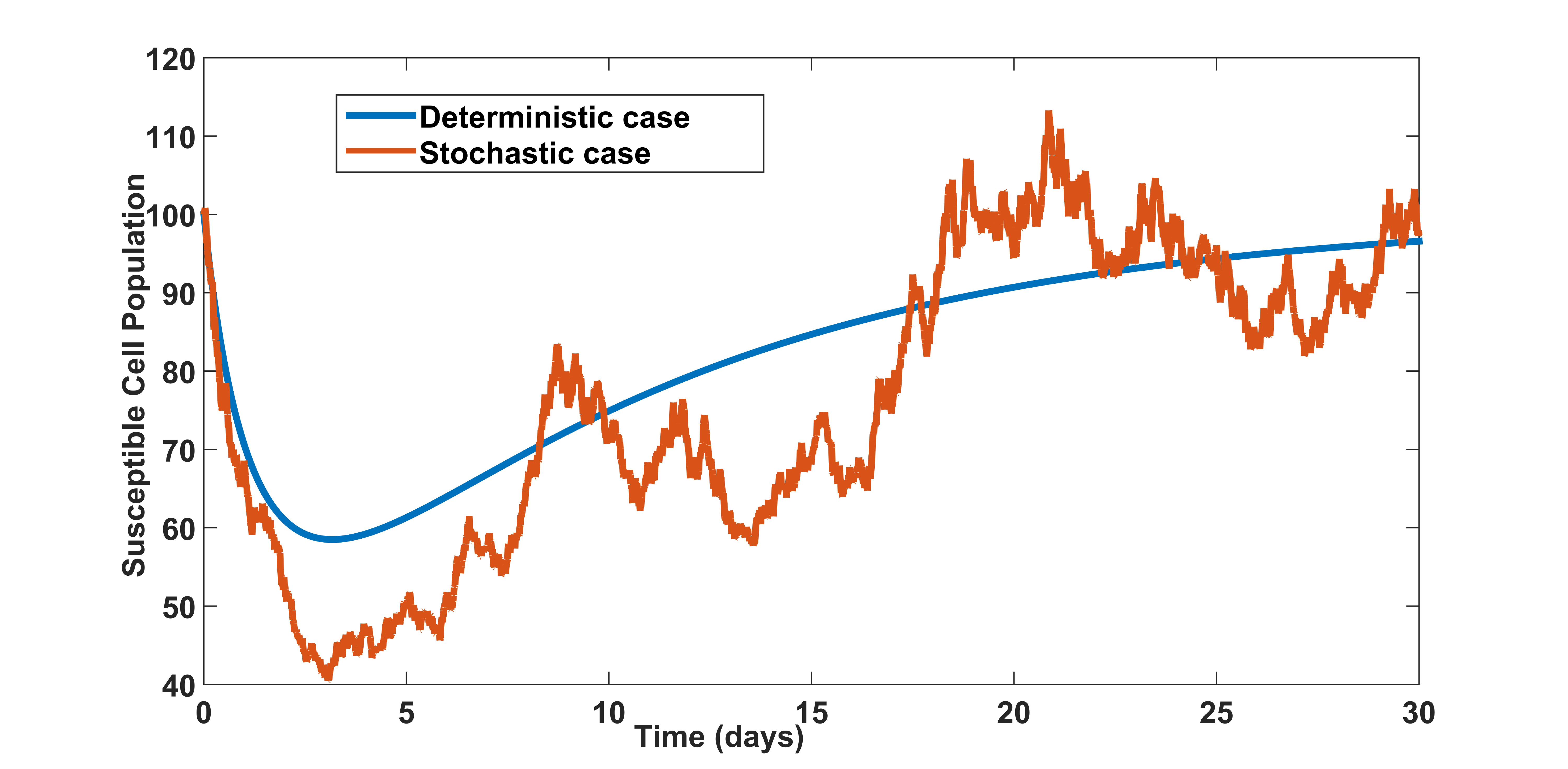}
\hspace{-.4cm}
\includegraphics[width=3.8in, height=3.3in, angle=0]{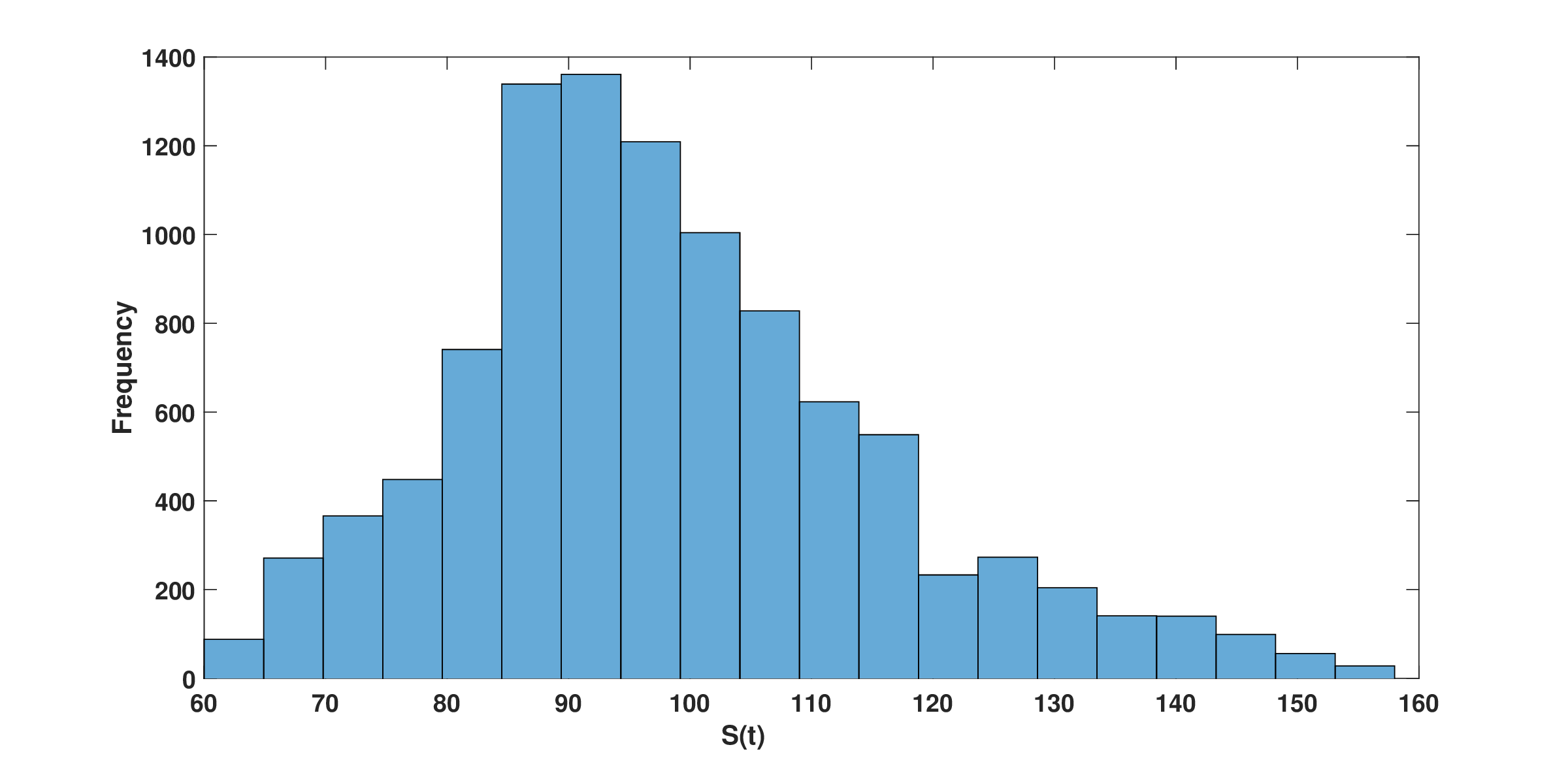}
\vspace{5mm}
\caption{Figure depicting the dynamics of susceptible population and distribution of susceptible population over time. We see that the variable $S(t)$ is distributed around the mean value of $\frac{\omega}{\mu}$ in case of stochastic model.  }
\label{fig2}
\end{center}
\end{figure}

\begin{figure}[hbt!]
\begin{center}
\includegraphics[width=6.8in, height=4.2in, angle=0]{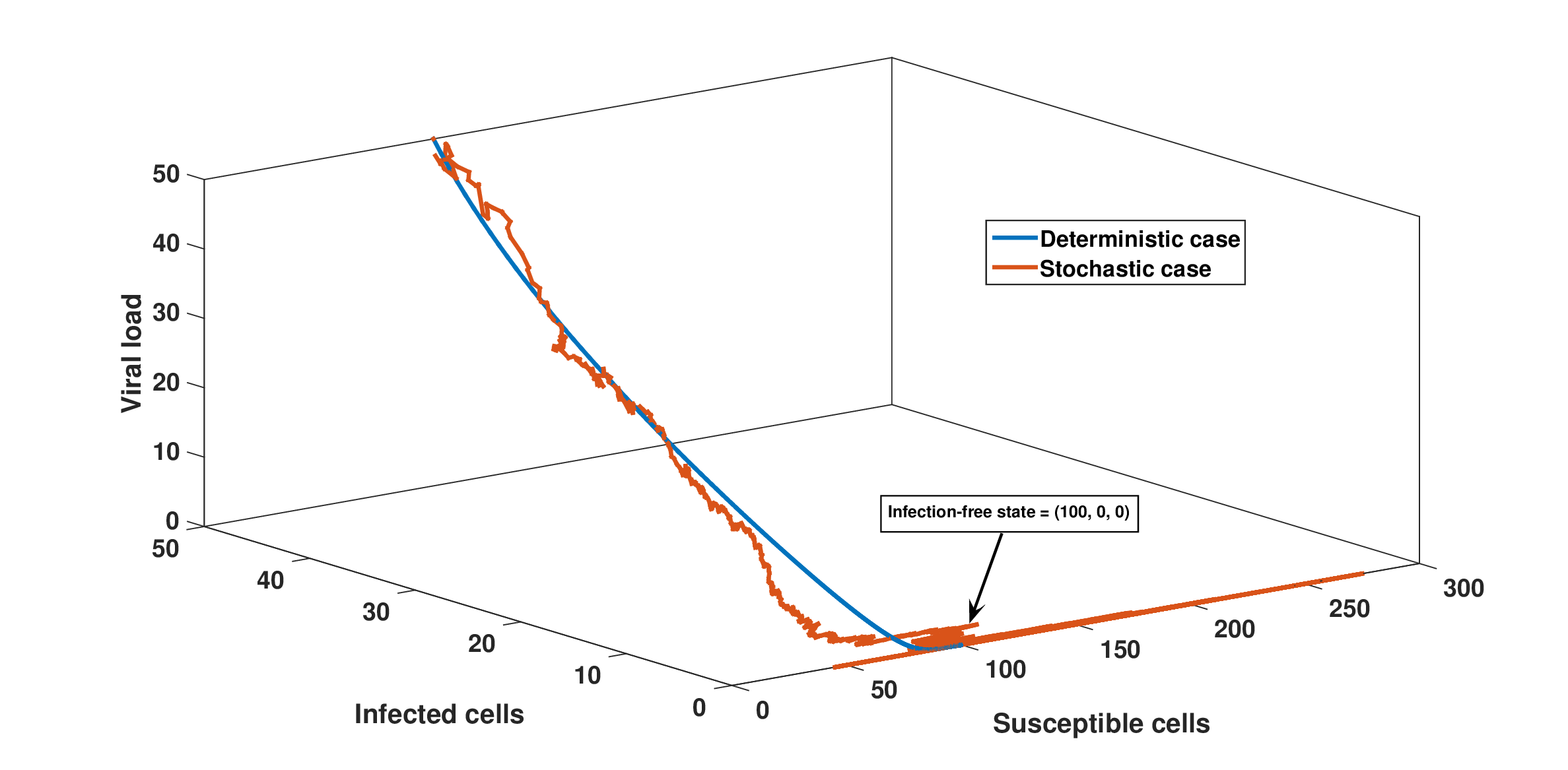}
\caption{Three-dimensional phase diagram of $S$, $I$ and $B$ for the case with $R_0 < 1$ and $\sigma_1 = \sigma_2 = 0.1$. We see that the population eventually stabilizes around the infection-free state $(100, 0, 0)$.}
\label{fig3}
\end{center}
\end{figure}

 \begin{example}
    \textnormal{ Taking the parameter values from Table \ref{table2} with $\sigma_1 = \sigma_2 = 0.1$,  we find that both the conditions A and B of Theorem \ref{thm38} are satisfied. For the deterministic case with the parameters values from Table \ref{table2}, the value of $R_0$ was found to be $0.15 < 1$. So in this case, according to the Theorem \ref{thm38}, the solution of the system tends to the infection-free equilibrium point at exponential rate, which means that the infected cell population and the viral load become extinct almost surely. This scenario is illustrated in Figure $\ref{fig1}$. Both deterministic and stochastic systems are plotted with initial value $(100, 100, 100)$ and we observe that initially the infected cell and viral load fluctuates in the presence of the environmental noise but eventually stabilizes around the infection-free state $(0, 0)$. } \\
    
    \textnormal{The susceptible cell population($S(t)$) for the deterministic case at the infection-free equilibrium point is given by $\frac{\omega}{\mu}$. With the parameter values from Table \ref{table2},  $\frac{\omega}{\mu} = 100$. In stochastic case, it is found that the variable $S(t)$ instead of being fixed at $\frac{\omega}{\mu} = 100$ is distributed around the mean value of $\frac{\omega}{\mu}$. This is illustrated in Figure $\ref{fig2}$. Three-dimensional phase diagram for this case is shown in Figure \ref{fig3} where we can see that the solutions eventually tend to $S$ axis. }
 \end{example}


\vspace{.4cm}

\begin{figure}[hbt!]
\begin{center}
\includegraphics[width=3.8in, height=3.3in, angle=0]{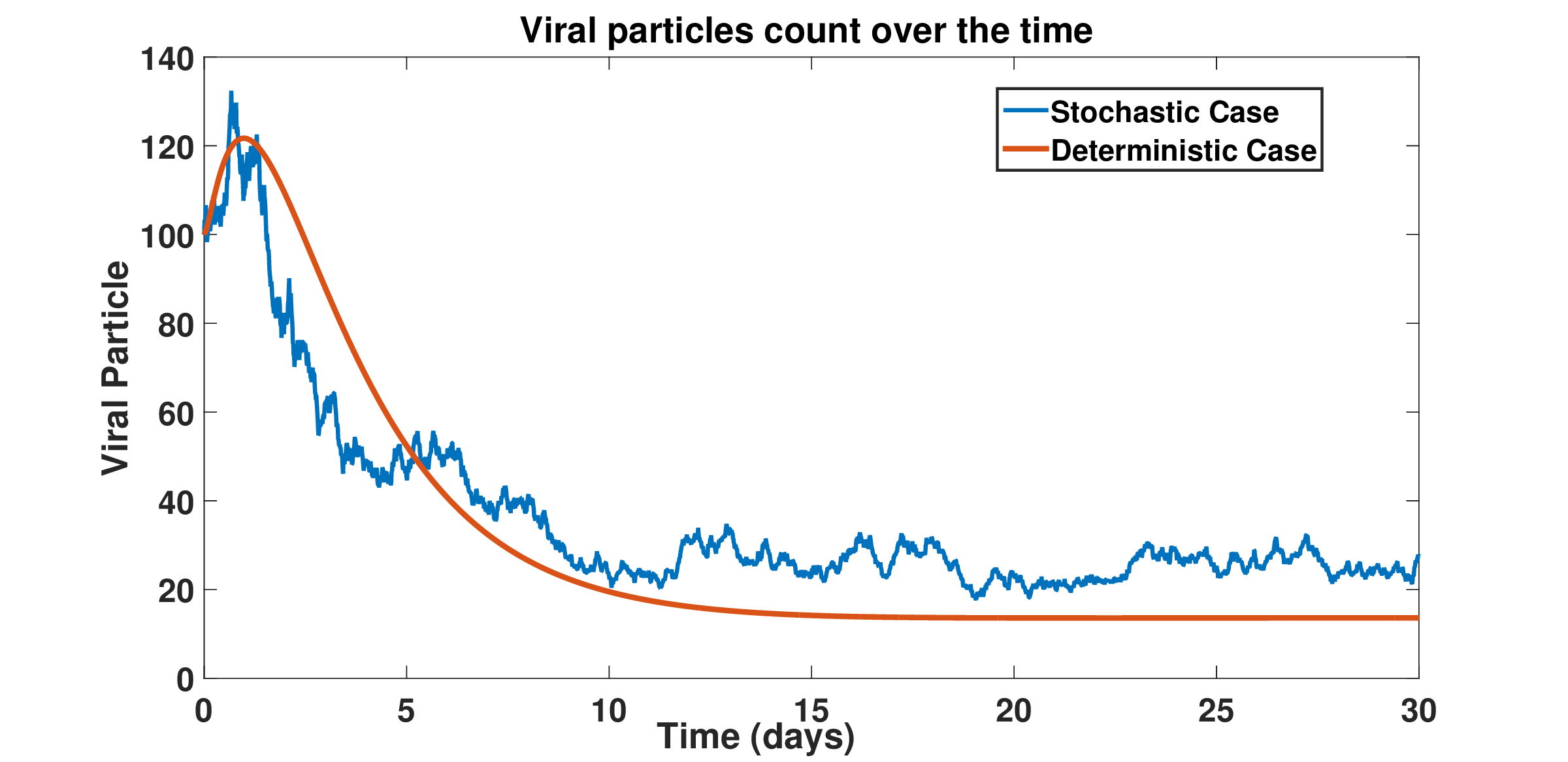}
\hspace{-.4cm}
\includegraphics[width=3.8in, height=3.3in, angle=0]{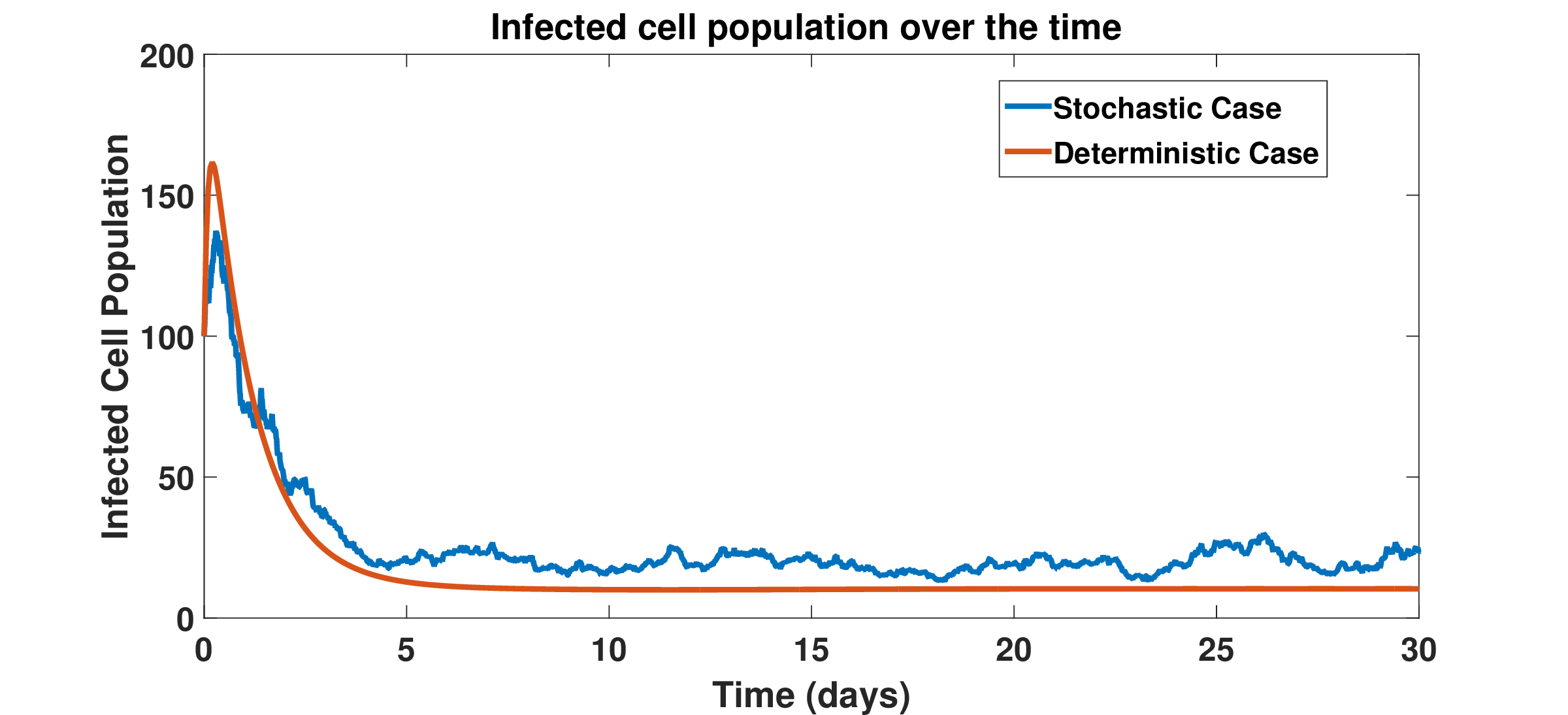}
\vspace{5mm}
\caption{Figure depicting the dynamics of  viral load population and infected cells for deterministic and stochastic cases for parameter values with $R_0 > 1$ and $\sigma_1 = \sigma_2 = 0.1$.}
\label{fig4}
\end{center}
\end{figure}

\begin{figure}[hbt!]
\begin{center}
\includegraphics[width=6.8in, height=4.2in, angle=0]{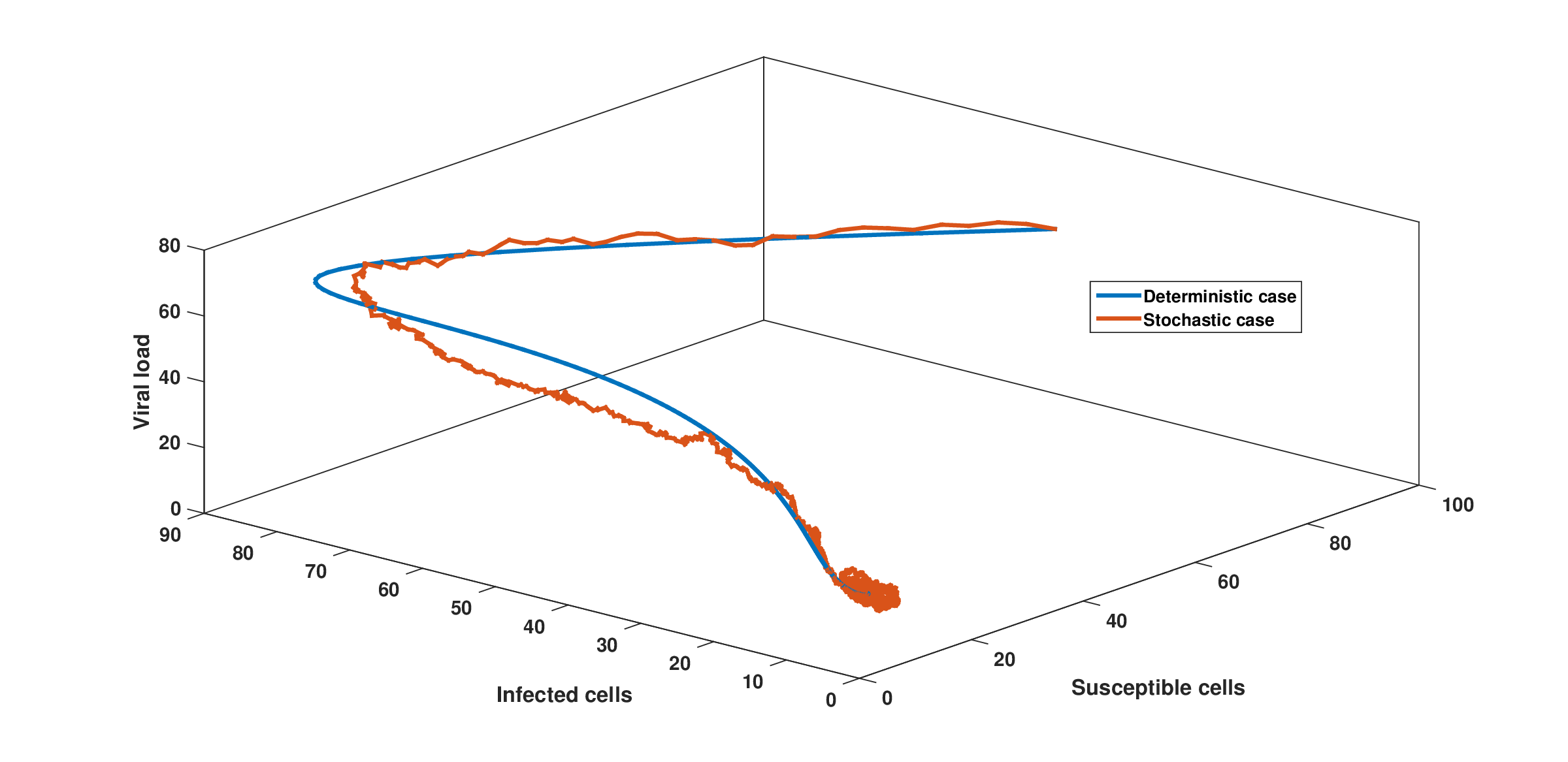}
\caption{Figure depicting susceptible, infected and viral load population over time. We see that the solution of the system eventually tend to a certain point in a three-dimensional space.}
\label{fig5}
\end{center}
\end{figure}

\vspace{.5cm}
\begin{example}
\textnormal{Figures $(\ref{fig1}, \ref{fig2}, \ref{fig3})$ are plotted with parameter values such that both the conditions of the Theorem \ref{thm38} are satisfied and in this case $R_0$ was less than unity. Now we want to slightly change the values of the parameters and see behaviour of the system. The parameter values in this case are taken as follows.}
$$\mu = \mu_1 = 0.1, \hspace{.15cm}\alpha = 0.24, \hspace{.15cm} \beta = 0.05,\hspace{.15cm} p=0.795, \hspace{.15cm}q=0.28, \hspace{.15cm}\beta =0.05$$
 \textnormal{and the initial value as $$(S_0, I_0, B_0) = (100, 100, 100)$$
 For these choice of parameter values, $R_0 > 1$ and both the conditions of Theorem \ref{thm38} are violated. This case is illustrated in Figure $\ref{fig4}$.  We simulate the equations $(\ref{equ4})-(\ref{equ6})$ with and without randomness and plot the infected and viral population over time. The value of noise intensities were taken as $\sigma_1 = \sigma_2 = 0.1$. In this case, we find that the infected cell and viral load count decreases over the time but does not die out.  In addition, we can see from the three-dimensional phase diagram (Figure \ref{fig5}) that the solution of the system tend to a certain point in a three-dimensional space. This case suggests that the infection will remain or persist with $R_0 > 1$ and $\sigma_1 = \sigma_2 = 0.1$.}  
 \label{example2}
 \end{example}
 
\begin{figure}[hbt!]
\begin{center}
\includegraphics[width=3.8in, height=3.15in, angle=0]{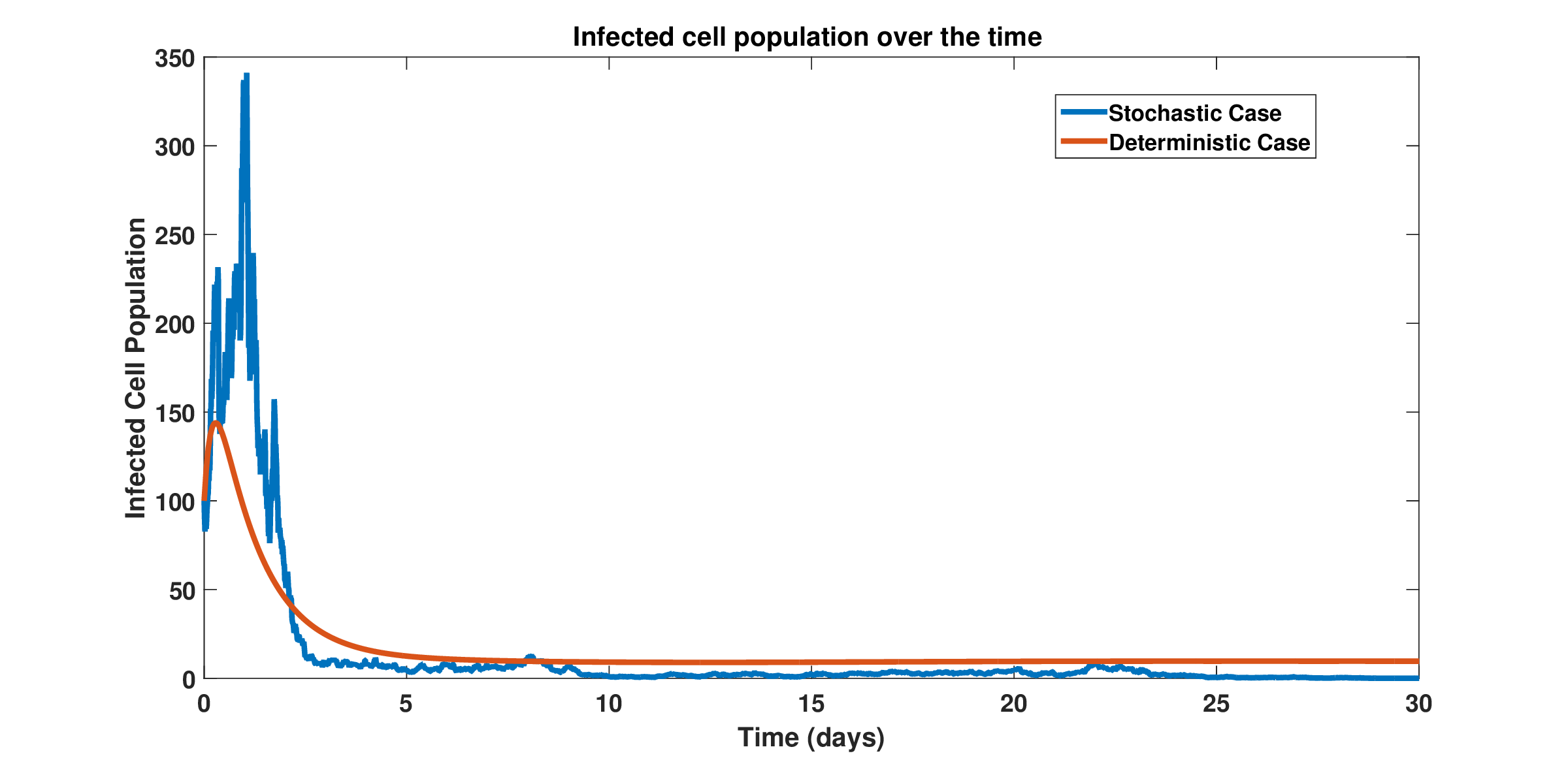}
\hspace{-.4cm}
\includegraphics[width=3.8in, height=3.15in, angle=0]{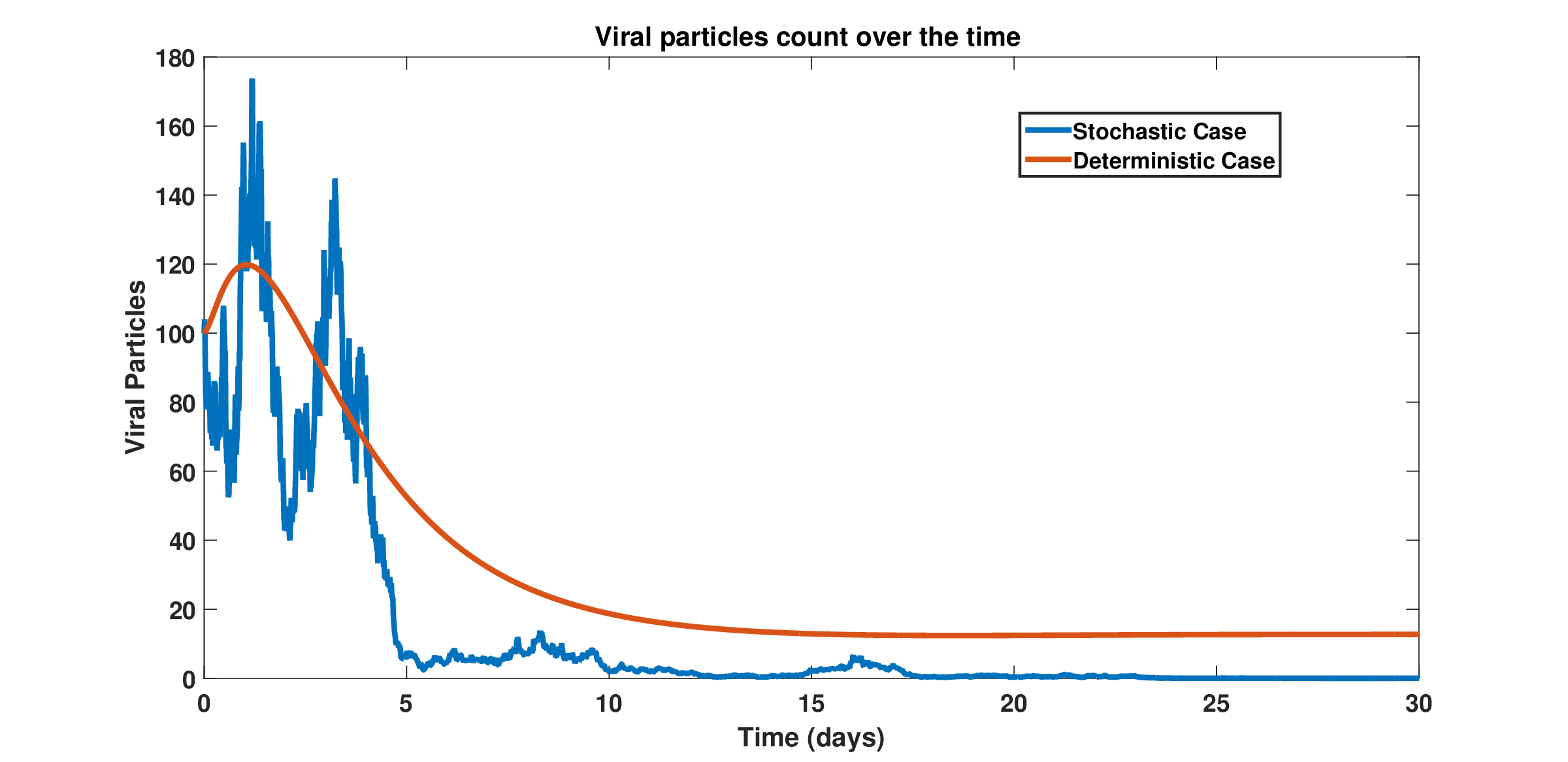}
\vspace{5mm}
\caption{Figure depicting the dynamics of  viral load population and infected for deterministic and stochastic cases for parameter values with $R_0 > 1$ and $\sigma_1 = 0.5$, $\sigma_2 = 0.8$. We observe that the infected cell and viral load vanishes under larger noise, although it survives for a long time in the deterministic
system.}
\label{fig6}
\end{center}
\end{figure}

 \vspace{.3cm}
 \begin{example}
\textnormal{Now what happens when noice intensities increases?. In this particular case,  we increased the noise intensities keeping parameters
and initial values same as that in Example \ref{example2}.  The values of $\sigma_1$ and $\sigma_2 $ were taken to be $0.5$ and $0.8$ respectively and with this the system of equations were solved. In this case, we observe (Figure \ref{fig6}) that the infected cell and viral load count fluctuates initially but eventually becomes zero. In other words, disease vanishes under larger noise, although it survive in the deterministic
system. This observation suggests that environmental noise plays a crucial role in controlling the spread of an epidemic. The disease can die out even when $R_0 > 1$. This means that larger noise can lead to the extinction of the disease, even though the disease is persistent in the deterministic model. }\\
\end{example}

\begin{figure}[hbt!]
\begin{center}
\includegraphics[width=3.75in, height=3.1in, angle=0]{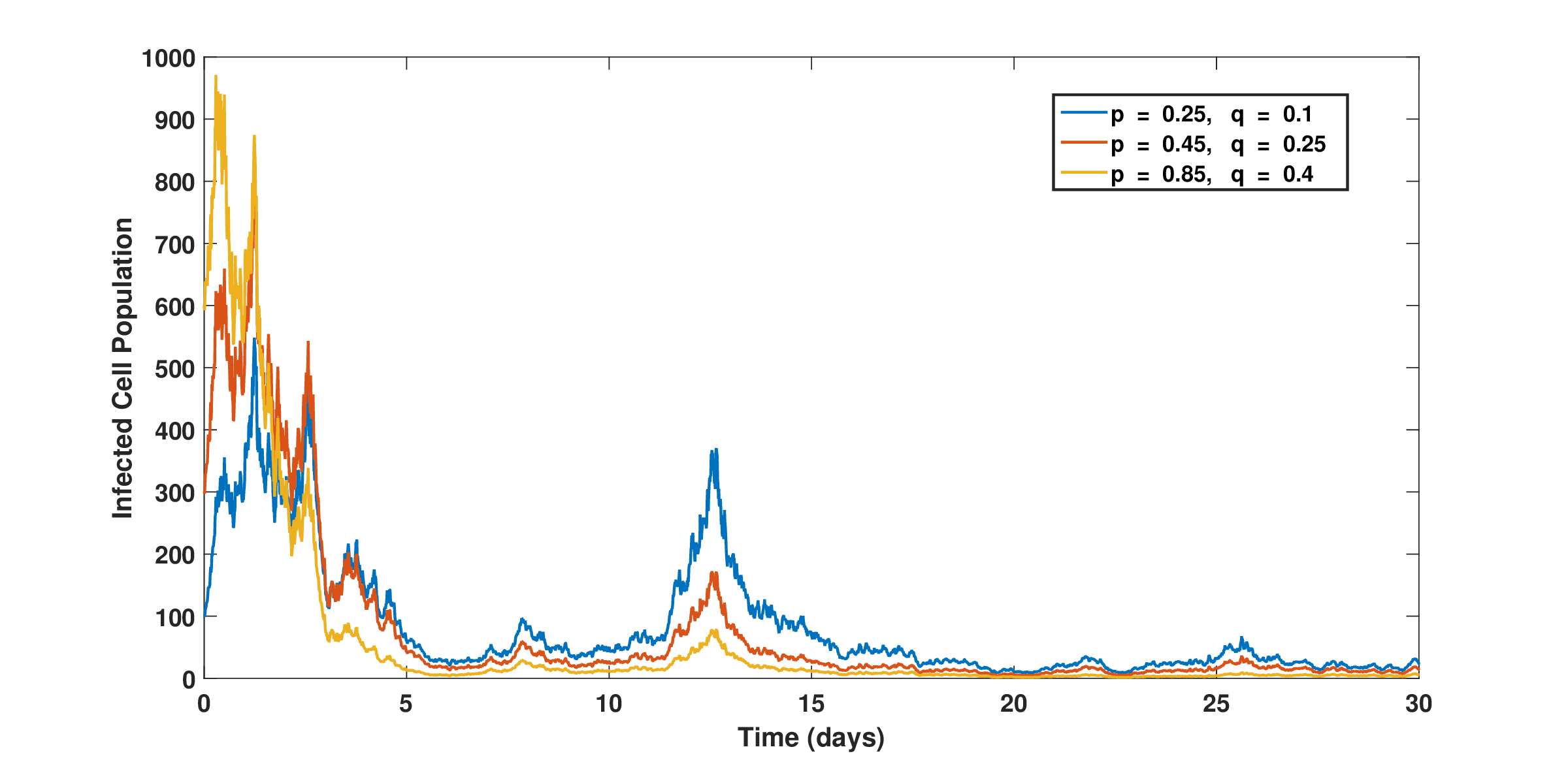}
\hspace{-.4cm}
\includegraphics[width=3.75in, height=3.1in, angle=0]{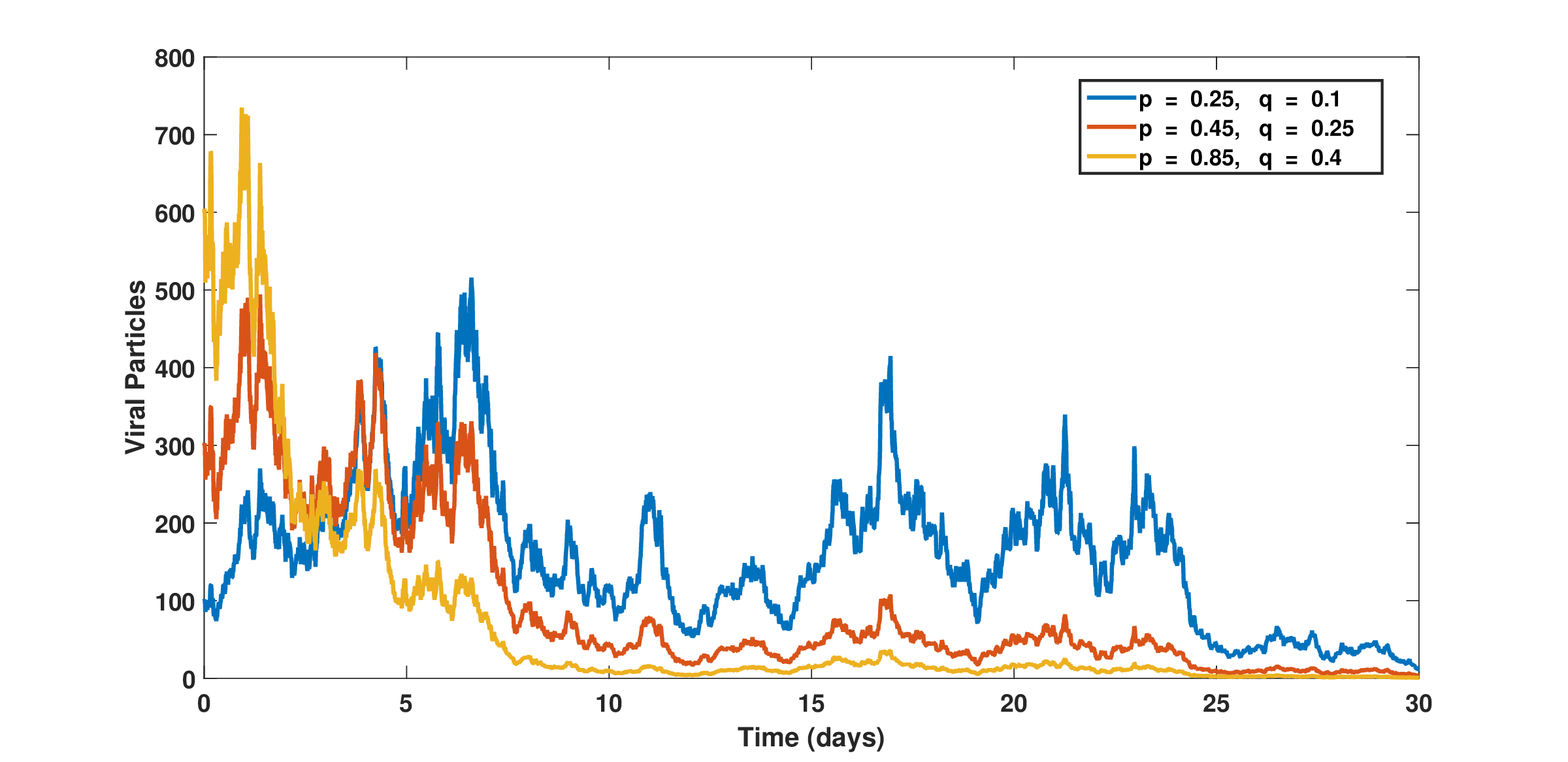}
\vspace{5mm}
\caption{Figure depicting the dynamics of infected and viral load population for different values of $p$ and $q$. We observe that higher values of $p$ and $q$ results in lower average values of the infected and viral load population.}
\label{fig7}
\end{center}
\end{figure}

\vspace{.2cm}
\begin{example}
\textnormal{The immune response to a virus or any foreign pathogens involves various parts of the immune system. Different types of immune cells are crucial in recognizing and attacking the virus. Cytotoxic T cells, natural killer cells, and macrophages activate other immune cells and directly or indirectly kill virus-infected cells \cite{azkur2020immune,burioni2021assessing}. In our model $(\ref{equ4})-(\ref{equ6})$,  the clearance rate of the infected cells and viral particles by the activation of the host immune system are represented by the parameters $p$ and $q$. In Figure $\ref{fig7}$,  we plot $I(t)$ and $B(t)$ by taking different values of $p$ and $q$ and the rest of the parameter values from Table \ref{table2}. In particular, we take three different cases, case $(i): p = 0.25, q = 0.1$, case $(ii): p = 0.45, q = 0.25$, and case $(iii): p = 0.8, q = 0.4$. To get a better perspective and global behaviour, we take different initial values of the population for each of these cases. For case $(i)$, $(S_0, I_0, B_0) = (200, 40, 100)$, for case $(ii)$, $(S_0, I_0, B_0) = (200, 300, 300)$, and for case $(iii)$, $(S_0, I_0, B_0) = (200, 600, 600)$. We observe that the average value of $B(t)$ in case (iii) is significantly lower compared to the other two cases. A similar pattern can be seen for I(t). This suggests that a stronger immune response leads to a more efficient clearance of infected cells and viral particles. Consequently, individuals with a robust immune system are likely to recover quickly from the infection.}
\end{example}

\vspace{1cm}
\section{Optimal Control Studies} \vspace{.2cm}
\subsection{ \textbf{Stochastic Optimal Control Problem Formulation} }
 
\vspace{.2cm}
In mathematical modelling studies optimal control problem focuses on determining the best control strategy for a dynamical system over a specified time period to optimize an objective function \cite{ross2009primer}. Optimal control theory has been extensively applied to epidemic models to design strategies for minimizing the total infectious burden through various control policies. Application of optimal control on stochastic models can be found in \cite{liu2021stochastic, witbooi2023stability, zhang2020sufficient, ishikawa2012stochastic}. In this section, we formulate stochastic optimal control problem by adding control variables in the system $(\ref{equ4})-(\ref{equ6})$.  The details of the control measures are as follows. The deterministic version of this has been already studied in \cite{chhetri2021within}.
\begin{enumerate}
		\item \textbf{Immunomodulating Drugs:} The\textcolor{white}{\lq\lq}first line of defence (also known as innate response) involves immune cells such as NK cells, macrophages, neutrophils,  and basophils. These immune cells primarily work to contain the pathogen. The acquired immune response, which leads to the activation of T and B lymphocytes, targets the specific antigen \cite{delves2020innate}. We use the control variable $u_{11}(t)$ to denote the interventions that stimulate the innate response and reduce number of infected cells. The interventions that enhances the B and T cells response and remove the viral particles is denoted by $u_{12}(t)$.\\
		
	\item \textbf{Antiviral Drugs:} The primary function of antiviral drugs is to target the viral replication and reduce the multiplication rate. These drugs, when administered after the onset of symptoms in a patient, can reduce the risk of infection for others by reducing viral shedding in the patient's respiratory secretions \cite{mitja2020use}. The control variable $u_2(t)$ denotes the antiviral effect of the drug Aribidol. Other pharmaceutical agents like Remdesivir \cite{52}, Lopinavir/ritonavir \cite{62},  Favipiravir \cite{82} are also believed to prevent replication of virus to some extend. 
\end{enumerate}

With this, we define the set of all admissible controls as follows. \\

	$U = \left\{(u_{11}(t),u_{12}(t),u_{2}(t)) : u_{11}(t) \in [0,u_{11} max] , u_{12}(t) \in [0,u_{12} max] , u_{2}(t) \in [0,u_{2} max] ,t \in [0,T] \right\}$\\

We assume that the control variables are measurable and bounded functions. On this basis, we now propose the optimal control problem with the aim of reducing the cost functional given by,
	\begin{equation}
		J = E\bigg[\int_{0}^{T} (I(t)+B(t)+A_{1}(u_{11}^2(t) + u_{12}^2(t)) + A_{2}u_{2}^2(t)) dt\bigg]   \label{obj}
	\end{equation} 
with the variables satisfying the system 
\begin{eqnarray}
   	{dS}& =&  \bigg(\omega \ - \beta SB  - \mu S\bigg)dt - \sigma_1 S dW_1(t)  \label{op1} \\
   	dI &=& \bigg(\beta SB - pI - u_{11}(t)I - \mu I \bigg) dt- \sigma_1 I dW_1(t)  \label{op2}\\ 
   	dB &=&  \bigg((\alpha - u_2(t)) I -qB - u_{12}(t)B - \mu_{1} B\bigg) dt - \sigma_2 BdW_2(t)\label{op3}
   \end{eqnarray}

The cost functional (\ref{obj}) contains the population of infected cells and the population of viral load as well as the terms related to the potential side effects or toxicity of the drugs. The term $A_{1}(u_{11}^2(t) + u_{12}^2(t))$ penalizes the excessive use of the immunomodulatory drug and balances the benefits of the drug against its potential side effects. A higher value of $A_1$ means a  higher penalty for using the drug and therefore caution should be exercised in administering the drug to avoid the side effects or toxicity \cite{joshi2002optimal}. Similarly, $A_{2} u_{2}^2(t)$ weighs the benefits of the antiviral drug against its potential side effects. Increase in the control value does not necessarily lead to a decrease in the number of infected cells and viral load. In this particular scenario, an increase in the control value actually corresponds to an increase in the cost or potential side effects of the drugs, as shown by the nonlinear control variables in the definition of the cost functional $(\ref{obj})$. The goal here is to minimize the expected infected cells and viral\textcolor{white}{\rq\rq}load in the body by finding the optimal drug with the least side effects.\\

\noindent
For simplicity we define $x(t) = [x_1, x_2, x_3]= [S, I, B]$ and with this we rewrite system $(\ref{op1})$ - $(\ref{op3})$ as,

\begin{eqnarray}
    dx(t) = f(x, u)dt + g(x)dW(t)
\end{eqnarray}
\noindent
where
$$f(x, u) = [f_1, f_2, f_3]$$
$$g(x) = [g_1, g_2, g_3]$$
$$W = [W_1, W_1, W_2]$$
and
\begin{equation*}
    \begin{split}
f_1(x, u) &= \omega - \beta x_1 x_3 - \mu x_1\\
f_2(x, u) &= \beta x_1 x_3 -px_2 - u_{11}(t)x_2 - \mu x_2\\
f_3(x, u) &= (\alpha - u_2(t))x_2 - qx_3 - u_{12}(t)x_3 - \mu_1 x_3\\
g_1(x) &= - \sigma_1 x_1\\
g_2(x) &= - \sigma_1 x_2\\
g_1(x) &= - \sigma_2 x_3
    \end{split}
\end{equation*}

Before we proceed to the calculation of optimal control solutions that minimize the cost functional, we will make a remark about the existence of optimal controls. The admissible control set $U$ defined above is compact. From Theorem \ref{thm1}, we know that the SDE $(\ref{op1})$ - $(\ref{op3})$ with each bounded control variables has a unique solution that are non-negative for all $t\geq 0$. Also, the integrand of the cost functional $J$ is convex on $u = (u_{11}, u_{12}, u_{2})$. With these conditions, the existence of optimal control minimizing $J$ is guaranteed \cite{liu2021stochastic, prakash2023stochastic}.\\

\noindent		
We now define the Hamiltonian function and use stochastic maximum principle \cite{aastrom2012introduction, liu} to obtain the optimal control solution.
\begin{eqnarray}
    H(x, u, m, n)  =  <f(x, u), m> + L(x, u) + <g(x), n>
\end{eqnarray}
Here 
 $L(x, u) = I(t)+B(t)+A_{1}(u_{11}^2(t) + u_{12}^2(t)) + A_{2}u_{2}^2(t)$. $<. , .>$ denotes a Euclidean inner product. $m = (m_1, m_2, m_3)$, $n =(n_1, n_2, n_3)$ are called the adjoint vectors. From the stochastic Maximum principle  the following relation holds:
\begin{eqnarray}
    dx^*(t) = \frac{\partial H(x^*, u^*, m, n)}{\partial m}dt + g(x^*)dW(t) \label{k1} \\   \vspace{.4cm}
    dm(t) = - \frac{\partial H(x^*, u^*, m, n)}{\partial x}dt + n(t)dW(t) \label{k2} \\ \vspace{.4cm}
    H(x^*, u^*, m, n) = \text{max} \hspace{.2cm} H(x^*, u, m, n), u \in U
\end{eqnarray}
The boundary conditions of the equations $(\ref{k1})-(\ref{k2}) $ are given by,
$$x^*(0) = [x_1^*(0), x_2^*(0), x_3^*(0)] = [S_0^*, I_0^*, B_0^*]$$
$$m(T) = [0, 0, 0]$$

To find the optimal $u^*$ that minimizes the cost functional $J$, we have to solve the boundary value problem $(\ref{k1})-(\ref{k2})$ with the above mentioned constraints. With the given initial conditions, we have to solve the state system $(\ref{k1})$ forward in time and the adjoint system $(\ref{k2})$ backward in time with $m(T) = [0, 0, 0]$. We must therefore solve a non-linear stochastic forward and backward differential equation.\\

\noindent
Using $(\ref{k2})$ and differentiating the Hamiltonian, we arrive at the following adjoint equations.

\begin{equation}
\begin{aligned}
	dm_1(t) &= \bigg[m_{1}(\beta x_3+\mu)- m_{2} \beta x_3 + \sigma_1 n_1 \bigg] dt + n_1 dW_1\\
	dm_2(t) &= \bigg[m_{2}(p+ u_{11}(t) + \mu)- m_{3} (\alpha - u_2(t)) + \sigma_1 n_2 - 1 \bigg] dt + n_2 dW_1\\
 dm_3(t) &= \bigg[m_{1}\beta x_1- m_{2} \beta x_1 + m_3(q+\mu_1 + u_{12}(t))+ \sigma_2 n_3  - 1\bigg] dt + n_3 dW_2\\
	\end{aligned}
	\end{equation}
\noindent
The boundary conditions being
$ m _{1} (T) = 0, \  m _{2} (T) = 0, \  m _{3} (T) = 0. $ Now, to obtain the optimal controls, we will use the Hamiltonian minimization condition 
\begin{eqnarray*}
    \frac{\partial H}{\partial u^*} = 0
\end{eqnarray*}
Differentiating the Hamiltonian and solving the equations, we obtain the optimal controls as 
\begin{eqnarray*}
	u_{11}^{*} &=& \min\bigg\{ \max\bigg\{\frac{m_{2}I}{2A_{1}},0 \bigg\}, u_{11}max\bigg\}\\
	u_{12}^{*} &= &\min\bigg\{ \max\bigg\{\frac{m_{3}B}{2A_{1}},0 \bigg\}, u_{12}max\bigg\}\\
	u_{2}^{*}& = &\min\bigg\{ \max\bigg\{\frac{m_{3}I}{2A_{2}},0 \bigg\}, u_{2}max\bigg\}
	\end{eqnarray*}

\section{\textbf{Numerical Simulation} } 

In this section, we run the numerical simulations to evaluate the effectiveness of the control measures. In doing so, we investigate how the dynamics of the system are influenced by the control variables. For the numerical simulation, we use the previously described Euler-Maruyama method. First, we solve the state system $(\ref{op1})-(\ref{op3})$ forward in time with the initial guess of the controls. Then, we solve the adjoint state system backward in time using the boundary conditions. After we have found the state variables and the adjoint variables, the values of the optimal control are updated and the process is repeated with the updated control variables. The parameter values are chosen so that the value of $R_0$ is greater than one (the previously discussed case) and $\sigma_1 = \sigma_2 = 0.05$. The values of $n_1$, $n_2$, and $n_3$ are taken to be $0.01$, $0.02$ and $0.03$ respectively. In this case, we expect the solution of the system to remain at a certain level for the no control case. The positive weights chosen are $A_{1}$ = 10, $A_{2}$ = 5. $A_{1}$ is chosen high compared to $A_2$ because the hazard ratio of INF or immunomodulators is comparatively higher than that of arbidol \cite{IFNhazard, yin2021antiviral, chhetri2022optimal}. The hazard ratio (HR) is important to determine the rate at which people treated with a drug suffer a particular complication per unit time compared to the control population \cite{spruance2004hazard}. We also vary the values of the weight constants and observe the changes in the population dynamics. \\

\begin{figure}[hbt!]
\begin{center}
\includegraphics[width=3.8in, height=3in, angle=0]{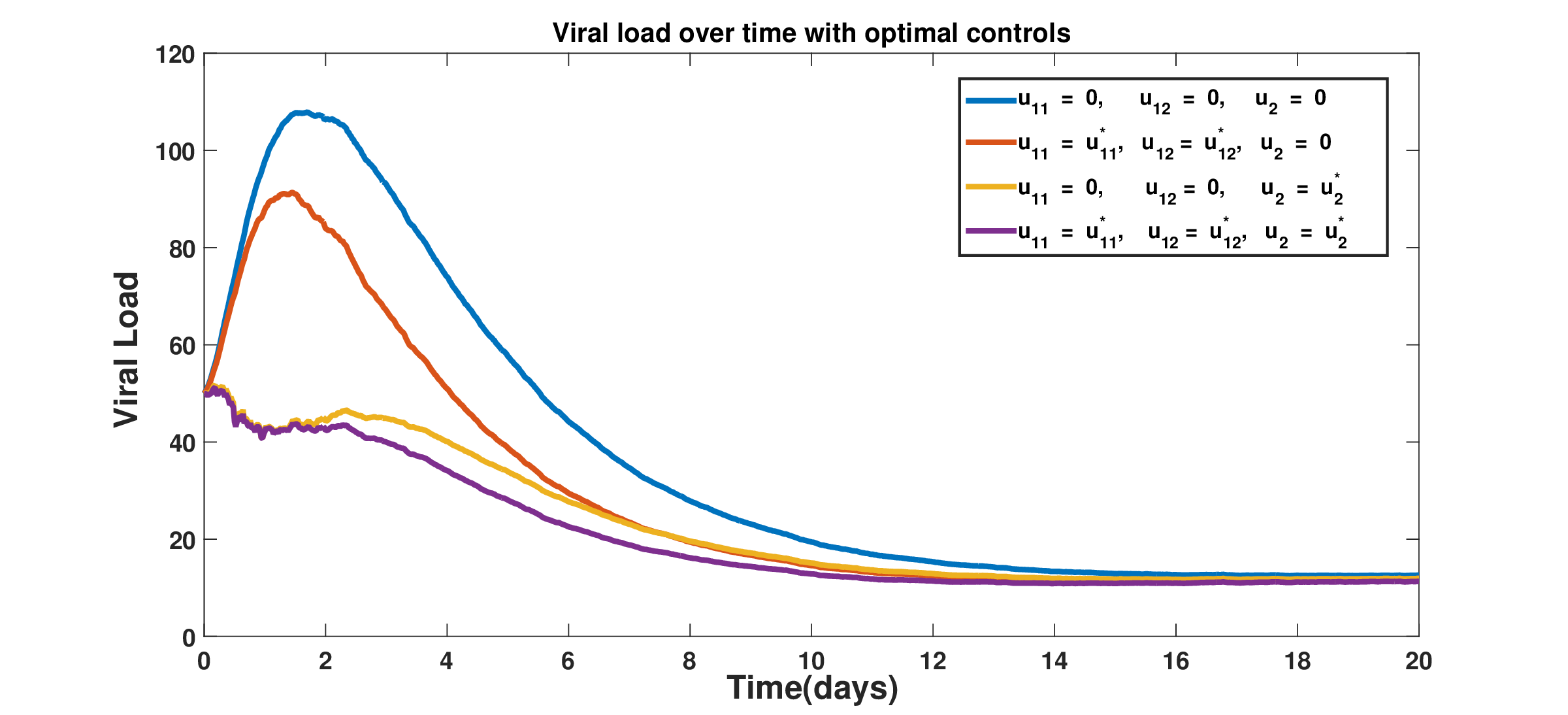}
\hspace{-.4cm}
\includegraphics[width=3.8in, height=3in, angle=0]{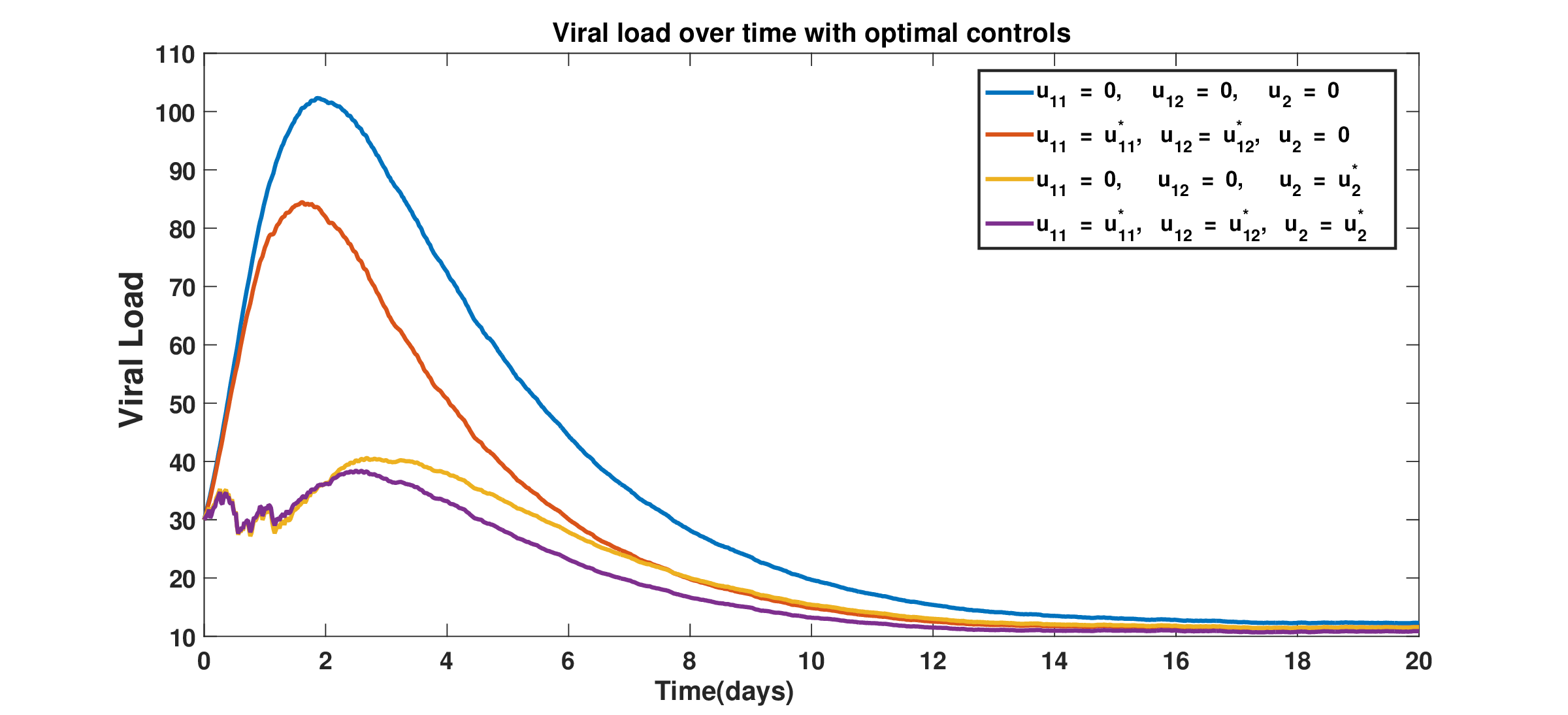}
\vspace{5mm}
 \caption{Figure depicting the dynamics of viral load ($B(t)$)  under the optimal controls $u_{11}^{*}, u_{12}^{*}$, $u_{2}^{*}.$ a) First figure is the case with $A_1 = 10, A_2 = 5$ and $(S_0, I_0, B_0)=(300, 80, 50)$. b) Second figure is the case with $A_1 = 10, A_2 = 5$ and $(S_0, I_0, B_0)=(200, 100, 30)$ }.
\label{fig66}
\end{center}
\end{figure}

\begin{example}
    
\textnormal{ In Figure \ref{fig66}, we have plotted the viral load population as a function of time for the individual and the combined control case. In the first Figure \ref{fig66} a) the initial value is assumed to be $(300, 80, 50)$ and in the second (Figure \ref{fig66} b)) it is $(200, 100, 30)$. In both cases, the value of the positive weight constant is assumed to be $A_1 = 10, A_2 = 5$. If no drugs are administered (case $u_{11} = u_{12} = u_2 = 0$, represented by the blue curve), there is initially a sharp increase in the viral load population, then the population decreases and stabilizes around the numerical value $10.5$. If only immunity-enhancing drugs are administered ($u_2 = 0$), a slight reduction in the peak value of the viral load can be observed compared to the case without control. The use of antiviral drugs $(u_2^*)$ shows better results compared to immunomodulators alone, but the combinations of these two drug interventions show much better results in minimising the viral load population in both cases. For the case with $u_{11} = u_{11}^*$, $u_{12} = u_{12}^*$, $u_2 = u_2^*$, we see that the viral load decreases and continues to decrease throughout the observation period, while for the other cases a slight increase in the number or a peak in the first few days are observed. There is not much difference observed in the behaviour of the viral load for the different initial conditions, except that in the second case (Figure  \ref{fig66} b) a slight increase in the count can be observed between $2^\text{nd}$ and $4^\text{th}$ day. But the overall pattern remains the same.}
\end{example}

\begin{figure}[hbt!]
\begin{center}
\includegraphics[width=3.8in, height=3in, angle=0]{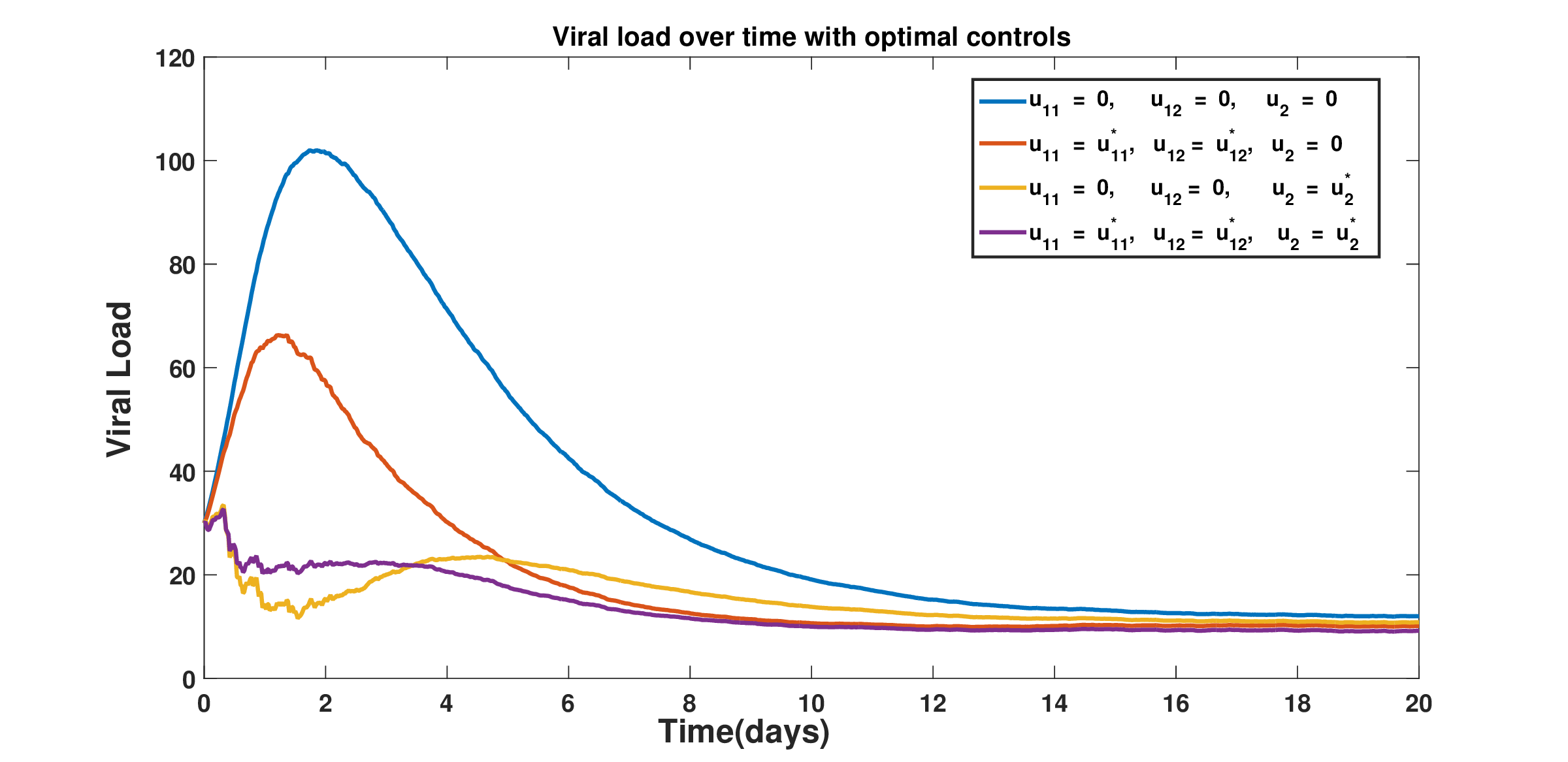}
\hspace{-.4cm}
\includegraphics[width=3.8in, height=3in, angle=0]{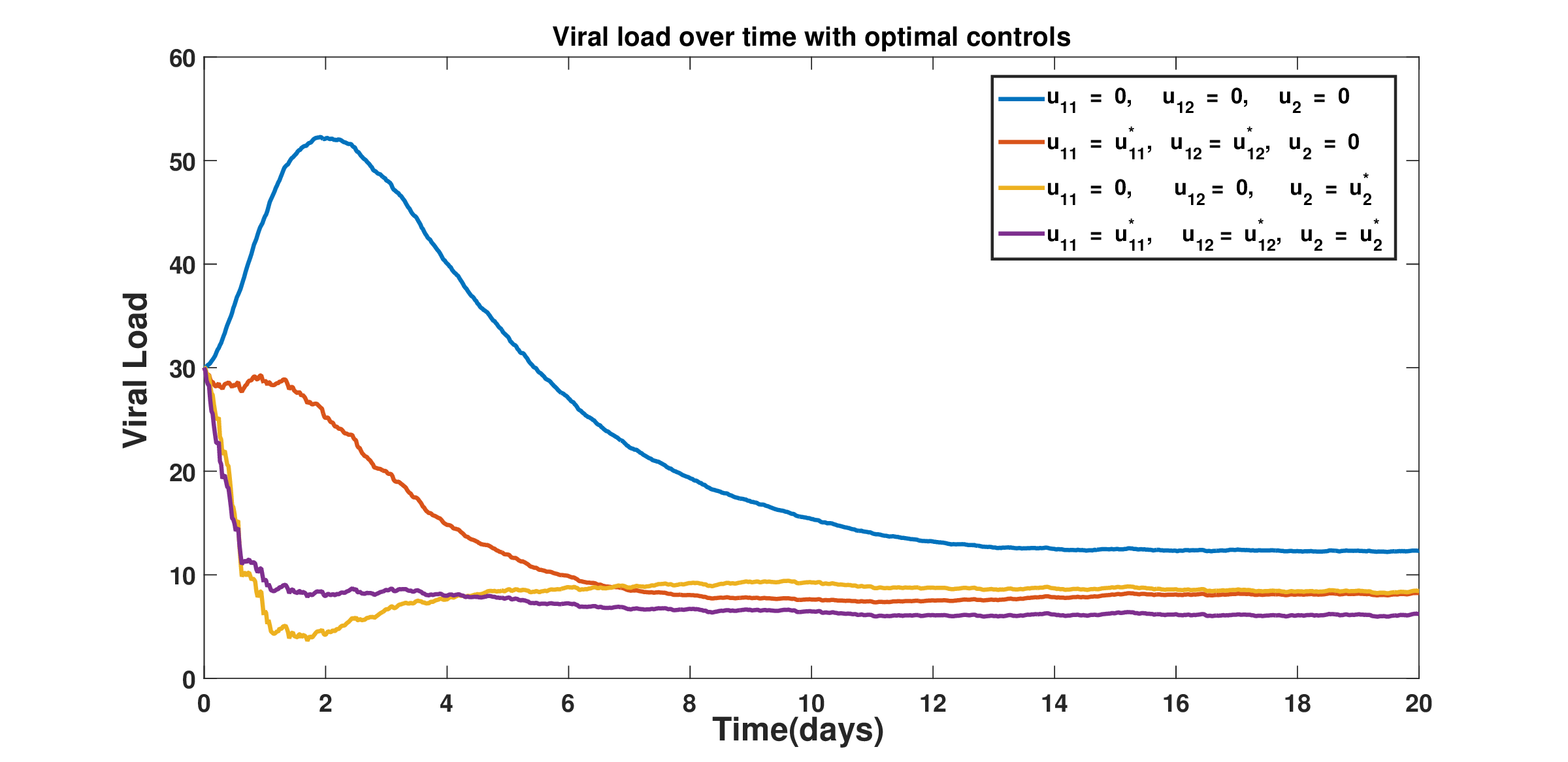}
\vspace{5mm}
 \caption{Figure depicting the dynamics of viral load ($B(t)$)  under the optimal controls $u_{11}^{*}, u_{12}^{*}$, $u_{2}^{*}.$ a) First figure is the case with $A_1 = 3, A_2 = 3$ and $(S_0, I_0, B_0)=(200, 100, 30)$. b) Second figure is the case with $A_1 = 1, A_2 = 0.8$ and $(S_0, I_0, B_0)=(200, 100, 30)$ }.
\label{figlkm}
\end{center}
\end{figure}

\begin{example}
\textnormal{ The definition of the optimal control solutions ($u_{11}^*, u_{12}^*, u_{2}^* $) shows that the optimal values of the control variables are in inverse proportion to the positive weight constants. The larger the values of $A_i$ for $i=1,2$,the smaller the values of the optimal controls. In Figure $\ref{figlkm}$, we have taken two different values for the weight constants, keeping the initial values of the state variables the same, and plotted the viral load population. In the first case, a), we have taken both $A_1$ and $A_2$ as $3$ (same values) and in the second case we have further reduced the values to $A_1 = 1$ and $A_2 = 0.8$. Compared to the previous Figure \ref{fig66}, we observe a slightly different behaviour of the viral population here. The drugs that prevent viral replication (denoted by $u_2^*$) initially prove effective, but due to the absence of immunomodulators or immune-enhancing drugs, their efficacy decreases and the viral load increases again. The viral load remains at a lower level with the combination of the controls.}\\

\textnormal{In both cases discussed above, we find that the optimal control strategy is to use a combination of both drugs to minimize the growth of the virus and thus accelerate the patient's recovery. This observation for the stochastic case is consistent with the results obtained for the deterministic case \cite{chhetri2021within}. Our observations are also consistent with the clinical results. Antiviral drugs such as Remdesivir, Arbidol are used in combination with Baricitinib, Tocilizumab etc. to fight different stages of the disease \cite{akinbolade2022combination}. In \cite{100}, it is shown that in a randomized study, people who received interferon-beta, ribavirin and lopinavir/ritonavir showed faster viral clearance and faster clinical improvement than a control group receiving only lopinavir/ritonavir.}
\end{example}

\section{Discussion and Conclusion} \vspace{.25cm}

This work is an extension of the study carried out in \cite{chhetri2021within}. In \cite{chhetri2021within}, the authors studied the dynamics of the COVID-19 pandemic a the virus cell level, including the optimal control problem using a deterministic approach. 
Biological systems are very complex and exhibit numerous interactions within the cell.
These systems are also subject to fluctuations in the environment. Environmental fluctuations play a crucial role in the development and progression of an epidemic \cite{bahar2004stochastic, li2015evolutionary}.  Temperature and rainfall variations have been found to lead to fluctuations in the dynamics of pathogenic fungi \cite{liu2015global}. The effects of absolute humidity on influenza virus transmission is investigated in \cite{shaman2009absolute}. It was found that absolute humidity significantly limits both transmission efficiency. The presence of even a tiny amount of white noise can suppress a potential population explosion \cite{mao2002environmental}. Therefore, it is important to study the effects of random fluctuations in the environment on population dynamics. In this study, a stochastic model is developed to describe the interactions between healthy cells, infected cells and viral particles. Natural mortality rates are assumed to be influenced by external noise, with the intensity of the noise quantified by intensity coefficients.\\

This work is divided into two parts. First, we establish the global positivity, boundedness and persistence of the model's solution and then examine the stability of the model. Our analysis shows that the number of infected cells and the viral load approach exponentially towards zero for higher values of the noise and $R_0 < 1$. The susceptible cell population on the other hand is distributed around the mean of its value at the infection-free equilibrium point. In scenarios where $R_0 > 1$ and the noise intensity is low, we observe that the number of infected cells and the viral load remain at a low level but do not vanish (indicating stabilization around the infected or disease equilibrium state). Interestingly, however, we find that at higher noise intensity, the number of infected cells and viral load drop to zero over time. This suggests that the disease can die out at higher noise, even when $R_0$ is greater than one (the disease persists in the deterministic setting). In \cite{elbaz2022modeling}, the authors have shown that the equilibrium state $E_0$ remains stable with increasing values of the noise parameter and the number of infected cells approaches zero. It is also found that the noise can even stabilize an unstable endemic deterministic system. Our results are consistent with these observations. \\

In the second part of the study, we formulate a stochastic optimal control problem and use the stochastic maximum principle \cite{witbooi2023stability, peng1990general} to find the optimal treatment rates that minimize the expected number of infected cells, viral particles, and costs or side effects caused by drug administration. The results of the optimal control studies show that both the antiviral drug that targets viral replication and the drug that enhances or modulates the immune response reduce the viral load when taken individually. Especially in the first few days of treatment, antiviral drugs that target viral replication seem to achieve better results than drugs that boost the innate immune response. However, it appears that these drugs achieve the best possible results when used in combination. The study suggests that the optimal control strategy could be to use a combination of the immunomodulating drug such as INF and the antiviral drug. This combination could potentially accelerate patient recovery.  Antiviral drugs such as Remdesivir, Arbidol are being used in combination with Baricitinib, Tocilizumab etc. to fight different stages of the disease \cite{akinbolade2022combination}. People who received interferon-beta, ribavirin and lopinavir/ritonavir are shown to have faster viral clearance and faster clinical improvement than a control group receiving only lopinavir/ritonavir \cite{100}. In \cite{kojima2022combination}, it is found that remdesivir along with immunomodulators can provide additional benefit in improving the respiratory status of COVID-19 patients. Combination therapy has been found to be effective not only in case of COVID-19 but also in several other infections. Some of the studies that discusses the effectiveness of combination therapy for HIV infection and hepatitis B virus can be found in \cite{joshi2002optimal, kirschner1997mathematical, yosyingyong2019global}. Our results are consistent with the clinical observations; however, comprehensive research is necessary, along with robust mathematical models that account for the potential adverse events of these drugs, to make well-informed decisions. \\

After the outbreak of COVID-19 at the end of 2020, several mathematical models were build and studied on different time scales. Most of these models have used deterministic equations and ignored the stochastic nature of systems. Our study shows that stochastic differential equations are effective in modeling virus dynamics and provide an alternative to traditional deterministic models. By replicating and improving the results of deterministic approaches, we highlight the added value of the stochastic model presented here. The inclusion of stochastic effects in models is crucial since most real-world problems are inherently non-deterministic, allowing for a more realistic representation of viral dynamics. In our study, we incorporated stochastic variations only in the within-host mortality parameters $\mu$ and $\mu_1$ of the COVID-19 model. Stochastic environmental variations can also be introduced and investigated for parameters such as the transmission rate and the birth rate, as these parameters are also strongly influenced by environmental noise.\\

\noindent
\textbf{Dedication} \\
The first author dedicates this research article to his loving father Late Purna Chhetri who still lives in his heart.\\

\noindent
\textbf{CRediT authorship contribution statement}\\
\noindent
BC conceived the idea and mathematical model, conducted the analysis, and numerical simulations.
BVRK contributed in conceptual discussion, organisation, planning, supervision, and corrections.\\

\noindent
 \textbf{Conflict of interest:} On behalf of both authors, the corresponding author states that there is no conflict of interest.\\
 
\noindent
\textbf{Funding:} Not Applicable\\

\noindent
 \textbf{Availability of data:} Not Applicable \\

\noindent
\textbf{Use of AI tools declaration}\\
The authors declare they have not used AI tools in the creation of this article.

  \bibstyle{plain}
	\bibliography{reference}
	
\end{document}